\documentclass[11pt]{article}
\usepackage{amsfonts, amsmath, amssymb, amsthm}

\usepackage[utf8]{inputenc}
\usepackage[colorlinks=true,citecolor=blue,urlcolor=blue,linkcolor=blue,bookmarksopen=true,unicode]{hyperref}
\usepackage{marginnote}

\title{Ranks of matrices with few distinct entries}
\author{Boris Bukh\thanks{Supported in part by U.S.\ taxpayers via NSF grant DMS-1301548.}}
\date{}

\makeatletter
\def\TODO{\@ifnextchar[{\TODO@with}{\marginnote{TODO}}}
\def\TODO@with[#1]{\marginnote{#1}}
\makeatother

\theoremstyle{plain}
\newtheorem{theorem}{Theorem}
\newtheorem{lemma}[theorem]{Lemma}
\newtheorem{corollary}[theorem]{Corollary}
\newtheorem{conjecture}[theorem]{Conjecture}
\newtheorem{proposition}[theorem]{Proposition}

\theoremstyle{remark}

\newcommand*{\abs}[1]{\lvert #1\rvert}                           
\newcommand*{\eqdef}{\stackrel{\text{\tiny{def}}}{=}}            
\newcommand*{\veps}{\varepsilon}                                 
\newcommand*{\Q}{\mathbb{Q}}                                     
\newcommand*{\R}{\mathbb{R}}                                     
\newcommand*{\Z}{\mathbb{Z}}                                     
\newcommand*{\F}{\mathbb{F}}                                     
\DeclareMathOperator{\rank}{rank}                                
\DeclareMathOperator{\vspan}{span}                               
\DeclareMathOperator{\Sec}{Sec}                                  
\DeclareMathOperator{\Gr}{Gr}                                    
\DeclareMathOperator{\Gal}{Gal}                                  
\DeclareMathOperator{\fchar}{char}                               
\DeclareMathOperator{\supp}{supp}                                
\DeclareMathOperator{\pred}{pred}                                
\newcommand*{\Proj}{\mathbb{P}}

\renewcommand{\theenumi}{\alph{enumi}}

\begin{document}

\maketitle

\begin{center}
\textit{In memory of a great teacher,\\Jirka Matou\v{s}ek}
\end{center}

\begin{abstract}
An $L$-matrix is a matrix whose off-diagonal entries belong to a set $L$, and whose
diagonal is zero. Let $N(r,L)$ be the maximum size of a square $L$-matrix of rank at most $r$.
Many applications of linear algebra in extremal combinatorics involve a bound on~$N(r,L)$.
We review some of these applications, and prove several new results on $N(r,L)$.
In particular, we classify the sets $L$ for which $N(r,L)$ is linear, and show that if $N(r,L)$ is superlinear
and $L\subset \Z$, then $N(r,L)$ is at least quadratic.

As a by-product of the work, we asymptotically determine the maximum
multiplicity of an eigenvalue $\lambda$ in an adjacency matrix of a digraph of a given size.

\end{abstract}

\tableofcontents

\section{Introduction}
\subsection{Motivation}
There are many applications of linear algebra to combinatorics that follow the same recipe.
They begin with $n$ objects of some kind, and a desire to bound~$n$. One then maps each of 
these objects to a pair $(v_i,u_i)\in V\times V^*$ where $V$ and $V^*$ are a vector space and 
its dual. The map is chosen so that the rank of the $n$-by-$n$ matrix $M=(u_i v_j)_{i,j}$ is 
large whenever $n$ is large. Since $\rank M\leq \dim V$, that yields a bound on~$n$.
In many of these applications $V$ is an inner product space, and $v_i=u_i$, but it is not always
the case.

The applications of this recipe include the proofs of the non-uniform Fisher inequality \cite{debruijn_erdos_prefisher,majumdar_nonuniform_fisher,isbell_nonuniform_fisher}, 
the Frankl--Wilson bound on $L$-intersecting families \cite{frankl_wilson}, Haemers' bound on the Shannon capacity of 
a graph \cite{haemers_shannon} and bounds on $s$-distance sets \cite{koornwinder_orig,blokhuis_fewdist}.
More applications can be found in the books by Babai--Frankl \cite{babai_frankl_book} and 
by Matou\v{s}ek \cite{miniatures}.

In all the applications named above, the matrices which arise are of a special form --- all diagonal
entries are equal, and the off-diagonal entries take on boundedly many distinct values. It is this
property that is used to bound their rank. The bounds on the ranks 
of such matrices are the subject of the present paper.

We define an \emph{$(L,\lambda)$-matrix} to be a square matrix whose diagonal entries are all equal to $\lambda$,
and each of whose off-diagonal entries is an element of the set~$L$. We shall mostly restrict
the study to $(L,0)$-matrices, which we call \emph{$L$-matrices} for simplicity. This incurs
only a minor loss of generality. Indeed, if $M$ is an $(L,\lambda)$-matrix and $J$ is the all-$1$ matrix, then $M-\lambda J$
is an $L'$-matrix for $L'=L-\lambda$, and the ranks of $M$ and $M-\lambda J$ differ by at most~$1$. The results
in this paper are too crude for this $\pm 1$ to matter. The advantage of the zero diagonal
is the dilation-invariance: if $M$ is an $L$-matrix, then $tM$ is an $tL$-matrix, for every scalar~$t$.

Suppose $L$ is a subset of some field, and $r$ is a natural number. We then define
\[
  N(r,L)=\max\{ n : \exists\ n\text{-by-}n\text{ }L\text{-matrix of rank }\leq r\}.
\]
Usually the underlying field will be clear from the context, but when confusion is possible
we shall write $L_{\mathbb{F}}$ to signify that $L$ is to be regarded as a subset of the field $\mathbb{F}$.

Throughout the rest of the paper, we shall only consider the case $0\not\in L$, since otherwise $N(r,L)=\infty$.
The case $\abs{L}=1$ is also easy, since then any $(L,\lambda)$-matrix is of the form $aI+bJ$, and so is
of rank at least~$n-1$. Furthermore, the determinant is $\det(aI+bJ)=a^{n-1}(a+bn)$ making it  
straightforward to tell when the rank is $n$ and when it is~$n-1$.

The first non-trivial case is $\abs{L}=2$. There is a natural correspondence between ranks
of $L$-matrices with $\abs{L}=2$ and multiplicities of eigenvalues 
of directed graphs. In our terminology, the adjacency matrices of directed graphs are just $\{0,1\}$-matrices.
If $M$ is an $\{0,1\}$-matrix with eigenvalue $\lambda$ of multiplicity $m_{\lambda}$, then
$M-\lambda I$ is a $(\{0,1\},-\lambda)$-matrix of rank $n-m_{\lambda}$. With the loss
of $\pm 1$ discussed three paragraphs above, that matrix is in turn equivalent to a 
$\{\lambda,\lambda+1\}$-matrix. Since every two-element set is a dilation of $\{\lambda,\lambda+1\}$
for a suitable $\lambda$, we can obtain any $L$-matrix with $\abs{L}=2$ this way, and the 
process is clearly reversible.

In view of the importance of adjacency matrices, we devote Subsection~\ref{subsec:ltwo} to
the case $\abs{L}=2$, in addition to the results for general~$L$ elsewhere in the paper. In the same subsection,
we also discuss eigenvalues of \emph{graphs}, which correspond to eigenvalues of \emph{symmetric} 
$\{0,1\}$-matrices.

\subsection{General remarks on upper bounds}
The results of this paper, which we will present in detail in Section~\ref{sec:results},
can be informally summarized as asserting that the order of magnitude of $N(r,L)$
is determined by a (possibly indirect) application of the following upper bound,
whenever $N(r,L)$ is not too large.
\begin{proposition}[Proof is in Section~\ref{sec:uppbounds}]\label{prop:basic}
Suppose $L$ is a $k$-element subset of some field, and $0\not\in L$.
Then the size of any $L$-matrix of rank $r$ is at most
\begin{equation}
\label{eq:basicupperbound}
  \frac{r^k}{k!}+O(r^{k-1}).
\end{equation}
\end{proposition}

Sometimes it is possible to combine \eqref{eq:basicupperbound} with a reduction modulo a prime. For example, 
if $M$ is a $\{1,3,8\}$-matrix over $\Q$, then $M\bmod 5$ is a $\{1,3\}$-matrix over the 
finite field $\mathbb{F}_5$. Since reduction modulo a prime may only decrease the rank, 
$N(r,\{1,3,8\}_{\Q})\leq  N(r,\{1,3\}_{\mathbb{F}_5})\leq \nobreak r^2/2+O(r)$.
As a special case of Theorem~\ref{thm:superlinear} we will show that in fact $N(r,\{1,3,8\}_{\Q})=\Theta(r^2)$. 
The proofs of our upper bounds on $N(r,L)$ can be viewed as a slightly
more sophisticated example of the same idea, where
we reduce modulo an ideal other than~$p\Z$. For example, one can obtain an upper 
bound on $N(r,\{1,\alpha\})$ by reducing modulo the ideal $(1-\alpha)\Z[\alpha]$ 
in the ring $\Z[\alpha]$.  However, we shall follow a more direct approach, inspired by \cite[Theorem 3.5(4)]{khot_thesis}, 
that avoids the language of ideals.

The simple combination of \eqref{eq:basicupperbound} with reductions modulo an ideal provides the only
known asymptotic upper bounds on $N(r,L)$ for a fixed~$L$. If we permit $L$ to vary, then it is possible
to prove \emph{relative bounds}. For example, if $L=\{1,1+\veps\}\subset \R$ for some small $\veps$, then
as $\veps\to 0$ the matrix tends to $J-I$, from which it is easy to deduce
that $N(r,\{1,1+\veps\})=N(r,\{1\})=r$ whenever $\veps<\veps_0(r)$.
More precise bounds for ranks of small perturbations of the identity matrix
have been established by Alon \cite{noga_perturbed}. In the same paper, he also gives numerous applications
of such bounds. A similar bound is also known in the special context of equiangular lines \cite[Theorem~3.6]{lemmens_seidel}.

The known asymptotic upper bounds that improve upon the upper bound \eqref{eq:basicupperbound} use the application-specific
structure of a matrix. I am aware of two applications where specialized arguments have been used. The first concerns $L$-intersecting 
families. A family $\mathcal{F}\subset 2^{[r]}$ of sets is $L$-intersecting if $\abs{A\cap B}\in L$ for any distinct sets $A,B\in \mathcal{F}$.
If $\mathcal{F}$ is also $k$-uniform, i.e., all sets are of size $k$, then the consideration of characteristic vectors
yields an $(L,k)$-matrix of rank at most $r$. The specific structure exploited
in the results about $L$-intersecting families concerns the intersection of more than two sets. For example, 
Deza, Erd\H{o}s and Frankl \cite{deza_erdos_frankl} proved a bound of the form $c_{k,L} r^{\abs{L}}$ on the cardinality of an $L$-intersecting
$k$-uniform family with the constant $c_{k,L}$ that is superior to the one in \eqref{eq:basicupperbound}. 
Frankl \cite{frankl_classification} determined the maximal size of an $L$-intersecting $k$-uniform family for $k\leq 7$
and for all possible $L$ with two exceptions. 

The second application, where \eqref{eq:basicupperbound} has been improved, involves spherical codes.
Delsarte, Goethals, Seidel \cite{delsarte_goethals_seidel} define a \emph{spherical $L$-code} to
be a set $\mathcal{C}$ of unit vectors in $\R^n$ such that $\langle v,u\rangle\in L$ for distinct $v,u\in \mathcal{C}$. The matrix
of inner products of vectors from $\mathcal{C}$ is an $(L,1)$-matrix and has the additional property of being positive definite.
This property was used for example in \cite{lemmens_seidel,neumaier,bukh_equiangular,keevash_sudakov_equiangular,balla_draxler_keevash_sudakov} to prove bounds on the number
of equiangular lines in $\R^n$ with a prescribed angle. It was also used in \cite{musin_twodistance}
to give bounds on spherical two-distance sets.

\section{Statement of results}\label{sec:results}
\subsection{Sets of linear growth} 
Our first result is a classification of $L$ for which $N(r,L)$
is as small as it can possibly be:
\begin{theorem}[Proof is in Section~\ref{sec:main}]\label{thm:minimal}
For a set $L=\{\alpha_1,\dotsc,\alpha_k\}$, the following three statements are equivalent:
\begin{enumerate}
\item $N(r-1,L)\geq r+1$ for some natural number $r$;
\item There exists a homogeneous polynomial $P$ with integer coefficients satisfying $P(1,\dotsc,1)=1$
  and $P(\alpha_1,\dotsc,\alpha_k)=0$;
\item There exists a constant $c>1$, which depends on $L$, such that
$N(r,L)\geq cr$ for all large~$r$.
\end{enumerate}
Furthermore, the limit $\lim_{r\to\infty} N(r,L)/r$ always exists (but might be infinite, see Theorem~\ref{thm:superlinear}).
\end{theorem}

In the special case $L=\{1,\alpha\}\subset\mathbb{C}$, the part (b) of the preceding theorem is equivalent to the assertion
that $1/(1-\alpha)$ is an algebraic integer.

In the case $\abs{L}=2$ the value of $\lim_{r\to\infty} N(r,L)/r$ is determined in Theorem~\ref{thm:quadratic} below.

In part (a), $r+1$ cannot be replaced by~$r$. For example, the $\{-1,1\}$-matrix $\left(\begin{smallmatrix}0&1&1\\1&0&-1\\1&1&0\end{smallmatrix}\right)$ is of rank~$2$, but $N(r,\{-1,+1\}_{\Q})\leq N(r,\{1\}_{\F_2})\leq r+1$.

Given the relations that a set $L=\{\alpha_1,\dotsc,\alpha_k\}$ satisfies, it is possible to verify if the condition in part (b)  holds.
Namely, let $I(L)$ be the homogeneous ideal in $\Z[x_1,\dotsc,x_k]$ consisting of the integer polynomials vanishing at $\alpha=(\alpha_1,\dotsc,\alpha_k)$.
If $I(L)$ is generated by $f_1,\dotsc,f_l$, then checking if the condition in part (b) holds amounts to checking
if $\gcd\bigl(f_1(1,\dotsc,1),\dotsc,f_l(1,\dotsc,1)\bigr)=1$. If instead of $I(L)$, we know only $I_{\Q}(L)$, which is the
homogeneous ideal in $\Q[x_1,\dotsc,x_k]$ of all the \emph{rational} polynomials vanishing at $\alpha$, then we can first compute $I(L)=I_{\Q}(L)\cap \Z[x_1,\dotsc,x_k]$
using the algorithm sketched in \cite{user}. 

\subsection{Sets of superlinear growth}
If the preceding theorem deals with those $L$ for which $N(r,L)=r+O(1)$, the next
one is about those for which $N(r,L)=O(r)$. To state it, let $L=\{\alpha_1,\dotsc,\alpha_k\}$ and
call $k$-tuple $(A_1,\dotsc,A_k)\in \Z^k$ a \emph{primitive linear relation} if
\begin{align*}
   A_1\alpha_1+\dotsb+A_k\alpha_k&=0,\\
   A_1+\dotsb+A_k&=1.
\end{align*}
Note that a primitive linear relation is a special case of polynomials appearing in part (b) of Theorem~\ref{thm:minimal}.
Namely, it is such a polynomial of degree~$1$. 

\begin{theorem}[Proof is in Section~\ref{sec:main}]\label{thm:superlinear}
For a set $L=\{\alpha_1,\dotsc,\alpha_k\}$, the following three statements are equivalent:
\begin{enumerate}
\item $N(r-1,L)>kr$ for some natural number $r$;
\item There exists a primitive linear relation on $L$;
\item $N(r,L)=\Omega(r^{3/2})$.
\end{enumerate}
Furthermore, if $\abs{L}\leq 3$, the exponent $3/2$ in (c) can be replaced by $5/3$. If $L\subset \Z$ or $\abs{L}=2$, the exponent $3/2$ can be replaced by $2$.
\end{theorem}
This result demonstrates several ways in which the function $N(r,L)$ is better behaved than the corresponding
extremal function for the problem of $L$-intersecting families. First, Frankl \cite{frankl_exponents} showed
that, for every rational number $s/d\geq 1$, there exists a set $L$ such that the maximum cardinality of an $L$-intersecting family
on $[r]$ is $\Theta(r^{s/d})$. Then, F\"uredi in \cite[Paragraph 9.3]{furedi_turan_survey}, extending an earlier work of Babai--Frankl \cite{babai_frankl_classification}, classified sets $L$ for which $L$-intersecting families have linear size, but the relevant condition on $L$ is computationally harder to verify than (b) above.
Finally, as shown by F\"uredi \cite{furedi_three_intersections} and Khot
\cite[Theorem~3.2]{khot_thesis}, for some $L$, the asymptotic size of largest $L$-intersecting families depends on the existence
of designs of a prescribed size and parameters. The problem of deciding whether a design with certain
parameters exists appears to be difficult, and many computational problems related to designs are known to be NP-hard \cite[p.~719]{crc_design_handbook}. 
It is likely that there is no similar obstruction to understanding $L$-matrices.

I conjecture that the exponent $3/2$ in Theorem~\ref{thm:superlinear} can be replaced by $2$ for all sets $L$ satisfying a primitive
linear relation.

\subsection{Eigenvalues of (di)graphs and the case \texorpdfstring{$\abs{L}=2$}{L=2}}\label{subsec:ltwo}
In this subsection we present a nearly complete determination of $N(r,L)$ for
two-element sets $L$. We also discuss the related problem of the maximum
multiplicity of a graph eigenvalue, which corresponds to the case
of symmetric matrices.

Here and throughout the paper, `multiplicity of an eigenvalue' refers to the geometric multiplicity. That is,
a matrix $M$ has eigenvalue $\lambda$ of multiplicity $m$ if the eigenspace associated to $\lambda$
is of dimension~$m$.

Let $F$ be a field, and $L\subset F$ be a two-element set. We
denote by $F_0$ the prime subfield of $F$, i.e., we define $F_0=\Q$
if $\fchar F=0$, and $F_0=\F_p$ if $\fchar F=p$.

As the value of $N(r,L)$ remains unchanged if we multiply elements of $L$
by a non-zero element of $F$, we may assume without loss of generality
that $L=\{1,\alpha\}$.

Let
\[
  E(n,\lambda)\eqdef \max \{ m : \exists\text{ an }n\text{-by-}n\ \{0,1\}\text{-matrix with eigenvalue }\lambda\text{ of multiplicity }m\}.
\]
If $M$ is a $\{1,\alpha\}$-matrix of rank $n-m$ then $\{0,1\}$-matrix
$(M+I-J)/(\alpha-1)$ has eigenvalue $1/(\alpha-1)$ of multiplicity $m-1$, $m$, or $m+1$.
Conversely, if $M$ is a $\{0,1\}$-matrix with eigenvalue $\lambda$ of multiplicity
$m$, then $M+\lambda(J-I)$ is a $\{\lambda,\lambda+1\}$-matrix of rank $n-m-1$, $n-m$, or $n-m+1$.
Hence,
\begin{equation}\label{eq:en}
\begin{aligned}
  E(n,\lambda)&=m      &&\implies &&N(n-m+1,\{\lambda,\lambda+1\})\geq n,\\
  N(r,\{1,\alpha\})&=n &&\implies &&E(n,1/(\alpha-1))\geq n-r-1. 
\end{aligned}
\end{equation}

If $M$ is a $\{0,1\}$-matrix, and $\alpha$ is an eigenvalue of multiplicity $m$,
and $p_\alpha$ is the minimal polynomial of $\alpha$ over $F_0$, then $p_{\alpha}^m$
divides the characteristic polynomial of~$M$. Hence, $m\leq n/\deg \alpha$, where
$\deg \alpha=\deg p_{\alpha}$ is the degree of $\alpha$ over~$F_0$. In view of the relation
between $E$ and $N$, we conclude that
\begin{equation}\label{eq:eigenupper}
\begin{aligned}
  E(n,\lambda)&\leq n/\deg \lambda,\\
  N(r,\{1,\alpha\})&\leq \left(1-\tfrac{1}{\deg \alpha}\right)^{-1}(r+1).
\end{aligned}
\end{equation}

The following result shows that the above bound is nearly tight. In
particular, for $\abs{L}=2$, it determines the limit as $r\to\infty$ of $N(r,L)/r$ in Theorem~\ref{thm:minimal}.
If $F$ is a field, we say that $\lambda$ is an \emph{algebraic integer} in $F$ if
$\lambda$ is a root of some monic polynomial $x^d+a_{d-1}x^{d-1}+\dotsb+a_0$ with integer
coefficients. The \emph{degree} of $\lambda$ 
is the least degree of such a polynomial. Note that if $\lambda$
is an eigenvalue of a $\{0,1\}$-matrix, then $\lambda$ is an algebraic integer.
\begin{theorem}[Proof is in Subsection~\ref{subsec:digraph}]\label{thm:quadratic}
Let $F$ be a field. Suppose $\alpha\in F$
is an element such that $\lambda=1/(1-\alpha)$ is an algebraic integer in~$F$.
Then
\begin{enumerate}
\item If $\deg \lambda=1$, then $E(n,\lambda)=n-\Theta(\sqrt{n})$ and $N(r,\{1,\alpha\})=\Theta(r^2)$;
\item If $\lambda$ has degree $d>1$, then
\begin{align*}
  \tfrac{d}{d-1}r-O(\sqrt{r}) &\leq N(r,\{1,\alpha\})\leq \tfrac{d}{d-1}(r+1),\\
  \tfrac{n}{d}-O(\sqrt{n})&\leq E(n,\lambda)\leq \tfrac{n}{d}.
\end{align*}
\end{enumerate}
\end{theorem}

\paragraph{Symmetric matrices and graph eigenvalues} 
We next discuss graph eigenvalues. We shall restrict our discussion to $\mathbb{C}$ as an ambient field.
The corresponding extremal functions are
\begin{align*}
  N_s(r,L)&\eqdef \max\{ n : \exists\text{ symmetric }n\text{-by-}n\text{ }L\text{-matrix of rank }\leq r\},\\
  E_s(n,\lambda)&\eqdef\max \{ m : \exists\text{ an }n\text{-vertex graph with eigenvalue }\lambda\text{ of multiplicity }m\}.
\end{align*}
The relations between $E$ and $N$ easily extend to $E_s$ and $N_s$, and we have
\begin{align*}
  E_s(n,\lambda)&=m      &&\implies &&N_s(n-m+1,\{\lambda,\lambda+1\})\geq n,\\
  N_s(r,\{1,\alpha\})&=n &&\implies &&E_s(n,1/(\alpha-1))\geq n-r-1.
\end{align*}

For a complex number $\lambda$ to be an eigenvalue of a symmetric integer matrix, $\lambda$ must be real.
Furthermore, as $\lambda$ is an algebraic integer and the Galois conjugates of $\lambda$ are 
eigenvalues of the same matrix, $\lambda$ must be in fact a totally real algebraic integer.

I conjecture that the extension of Theorem~\ref{thm:quadratic} to symmetric matrices holds for totally 
real algebraic integers.
\begin{conjecture}\label{conj:real}
Suppose $\lambda\in\mathbb{C}$ is a totally real algebraic integer of degree $d>1$, then
\[
  E_s(n,\lambda)\geq n/d-o(n).
\]
\end{conjecture}
We prove the conjecture for the degrees $d\leq 4$, and also for all `representable' $\lambda$.
We call a totally real algebraic integer $\lambda$ \emph{representable}
if there exists an integral symmetric matrix $M$ such that the map $\lambda\mapsto M$
is the isomorphism of algebras $\Z[\lambda]$ and $\Z[M]$. In other words,
the only eigenvalues of $M$ are $\lambda$ and its conjugates.

\begin{theorem}[Proof is in Subsections~\ref{subsec:digraph} and~\ref{subsec:graph}]\label{thm:reprconj}
If a totally real algebraic integer $\lambda$ is representable, then
Conjecture~\ref{conj:real} holds for $\lambda$.
\end{theorem}
\begin{theorem}[Corollary C in \cite{estes_guralnick}]
Every totally real algebraic integer of degree $d\leq 4$ is representable.
\end{theorem}

Estes and Guralnick \cite{estes_guralnick} conjectured that every
real algebraic integer is representable. I made the same conjecture
in a previous version of this paper. However, 
Dobrowolski \cite{dobrowolski_symmetric} disproved the conjecture.
Later, McKee \cite{mckee_six} constructed counterexamples of degree~$6$.

It is known that every totally real algebraic integer is an eigenvalue
of some symmetric matrix \cite{estes_symmetric,salez_tree}.

\subsection{Sets of arbitrary growth}
Whereas we do not have a complete classification of sets $L$
according to the growth rate of $N(r,L)$, we have two results
that restrict the possible growth rates.

Recall that a primitive linear relation on $L=\{\alpha_1,\dotsc,\alpha_k\}$ is an integer linear relation
of the form $\sum_i A_i \alpha_i=0$ with $\sum_i A_i=1$. Our first result
is that the growth rate of $N(r,L)$ is determined solely by the primitive linear relations
that $L$ satisfies. 
\begin{theorem}[Proof is in Section~\ref{sec:main}]\label{thm:lineareq}
Let $F$ be a field. Suppose $L,L'\subset F$ are sets of the same size $k$.
If $L$ and $L'$ satisfy the same set of primitive linear relations, then $N(r,L)\leq N(2kr,L')$, and similarly
$N(r,L')\leq N(2kr,L)$.
\end{theorem}

The second result is a generalization of the implication (a)$\implies$(b) from
Theorem~\ref{thm:minimal}.  Recall that a multivariate polynomial $P$ is said to vanish to order $m$
at a point $\alpha$ if all the monomials of degree at most $m-1$ in polynomial $P(x+\alpha)$ have zero coefficients. 

\begin{theorem}[Proof is in Section \ref{sec:uppbounds}]\label{thm:genupper}
Suppose $r,l,v$ are positive integers.
If $N(r-1,L)\geq \binom{r+l-1}{l}+\nobreak v$, then
there exists a $k$-variable homogeneous polynomial $P$ with integer coefficients
that satisfies the following
\begin{itemize}
\item $P(1,\dotsc,1)=1$;
\item For every univariate polynomial $f$ of degree at most $l$ with $f(0)=0$, polynomial 
$P$ vanishes at the point $\bigl(f(\alpha_1),\dotsc,f(\alpha_k)\bigr)$ to order at least $v$;
\item $\deg P\leq \binom{r+l-1}{l}+v$.
\end{itemize}
\end{theorem}

\section{Upper bounds}\label{sec:uppbounds}
The following lemma and its corollary, Proposition~\ref{prop:basic}, have been rediscovered several times
\cite{grolmusz_lowrank}, \cite{frankl_wilson} \cite[Lemma 2.3]{noga_perturbed}, and their origin is unclear.
The earliest references appear to be \cite{koornwinder_orig} and \cite{frankl_wilson}.
\begin{lemma}\label{lem:uppbound}
Suppose $M$ is a matrix over a field, $f$ is a univariate degree-$k$ polynomial, and $f[M]$ is the matrix
obtained from $M$ by applying $f$ to each entry. Then
\[
  \rank f[M]\leq \binom{\rank M+k}{k}.
\]
More generally, if $f$ contains only terms of degrees $d_1,\dotsc,d_t$, then the bound
is
\[
  \rank f[M]\leq \sum_i \binom{\rank M+d_i-1}{d_i}.
\]
\end{lemma}
\begin{proof}
In this proof, we write $vu$ for the vector that is the coordinate-wise product of 
vectors $v$ and~$u$, and $v^k$ for the coordinate-wise power of~$v$. 
Let $v_1,\dotsc,v_n$ be the columns of $M$. Put $r=\rank M$, and
assume without loss that $v_1,\dotsc,v_r$ span the column space of~$M$.
Let $V_d$ be the span of all the vectors that are of the form $v_1^{e_1}\dotsb v_r^{e_r}$
for some nonnegative exponents $e_1,\dotsc,e_r$ satisfying $e_1+\dotsb+e_r=d$.
Let $V$ be the span of $V_{d_1},\dotsc,V_{d_t}$.
Since $v_i=\sum_{j=1}^r \alpha_{ij} v_j$ for some scalars $\alpha_{ij}$,
it follows that $v_i^d=\sum_{j_1,\dotsc,j_d} \prod_{l=1}^d \alpha_{i j_l} v_{j_l}$.
In particular, $v_i^d\in V_d\subset V$. Hence each column of $f(M)$ 
lies in $V$. Since $\dim V\leq \sum_i \dim V_{d_i}=\sum \binom{r+d_i-1}{d_i}$, the result follows.
\end{proof}

The preceding implies an upper bound $N(r,L)\leq r^k/k!+O(r^{k-1})$ for every
$k$-element set $L$, i.e., Proposition~\ref{prop:basic}.
\begin{proof}[Proof of Proposition~\ref{prop:basic}]
We select a polynomial $f$
of degree $k$ that vanishes on $L$. Since $0\not\in L$, the matrix $f(M)$ is a non-zero
multiple of~$I$. The bound then follows from Lemma~\ref{lem:uppbound}.
\end{proof}

We next prove Theorem~\ref{thm:genupper}, which gives a necessary condition
for $N(r,L)$ to be large via the vanishing of a certain homogeneous polynomial.
For that we need a well-known lemma about the vanishing of the determinant function.

\begin{lemma}
Let $F$ be a field, and let $\det\colon M_n(F)\to F$ be the determinant. Regard
$\det$ as a polynomial in the~$n^2$ matrix entries.
Suppose $M$ is a matrix of rank at most $n-v$. Then $\det$ vanishes at $M$
to order at least $v$.
\end{lemma}
\begin{proof}
The proof is by induction on~$v$. 
The condition implies that every $(n-v+1)$-by-$(n-v+1)$ minor of $M$
vanishes. The partial derivative of $\det$ with respect to a matrix
entry is the cofactor of that entry. However, if $\fchar F\neq 0$, that is not
enough to complete the proof, as the vanishing of all the partial derivatives of orders up to 
$v$ is not equivalent to the vanishing of the polynomial to order~$v$. 

The rescue comes from the notion of a Hasse derivative. For a good
exposition of Hasse derivatives the reader might consult \cite[Section~2]{dvir_kopparty_saraf_sudan_multiplicities}.
Here, we recall only what we need. First, given a polynomial $P(x_1,\dotsc,x_k)$ and a multiindex $i\in \Z_+^k$, the Hasse
derivative $P^{(i)}$ of $P$ is defined as the coefficient of $z^i$ in $P(x+z)$, i.e., 
\[
  P(x+z)=\sum P^{(i)}(x)z^i.
\]

The induction proof then goes through in view of the following facts:
\begin{enumerate}
\renewcommand*{\theenumi}{H\arabic{enumi}}
\item The first Hasse derivatives are equal to the usual
first derivatives;
\item All Hasse derivatives of order at most $v-1$ vanish at a point if and only if the polynomial vanishes
to order $v$ at that point;
\item The Hasse derivatives satisfy the composition rule $(P^{(i)})^{(j)}(x)=c_{i,j}P^{(i+j)}(x)$,
where $c_{i,j}$ is a constant that we will not define here \cite[Proposition 4]{dvir_kopparty_saraf_sudan_multiplicities}.
\end{enumerate}

We use these facts to complete the proof. The base case $v=0$ is vacuous. 
Suppose $v\geq 1$. The induction hypothesis and (H1) tell us that all the first Hasse 
derivatives of $\det$ vanish to order at least~$v-1$. From (H2) and (H3), we then infer that
$\det$ vanishes to order at least~$v$.
\end{proof}
\begin{proof}[Proof of Theorem~\ref{thm:genupper}]
We first prove the case $l=1$. To that end, suppose
that $N(r-1,L)\geq r+v$, and $M_{r+v}$ is an $L$-matrix of size $r+v$
and rank $r-1$. 
We can regard $\det M_{r+v}$ as a homogeneous polynomial in $\alpha_1,\dotsc,\alpha_k$ of degree $r+v$,
say $\det M_{r+v}=P_{r+v}(\alpha_1,\dotsc,\alpha_k)$. By the preceding lemma
$P_{r+v}$ vanishes to order $v+1$ at $(\alpha_1,\dotsc,\alpha_k)$. Note that $P_{r+v}(1,\dotsc,1)=\det (J-I)=(-1)^{r+v-1} (r+v-1)$. 
Let $M_{r+v-1}$ be any principal $(r+v-1)$-by-$(r+v-1)$ submatrix of $M_{r+v}$. 
Similarly to the definition of $P_{r+v}$ we define homogeneous polynomial $P_{r+v-1}$ via 
$\det M_{r-v+1}=P_{r-v+1}(\alpha_1,\dotsc,\alpha_k)$. The polynomial vanishes to order
$v$ at $(\alpha_1,\dotsc,\alpha_k)$ and satisfies $P_r(1,\dotsc,1)=(-1)^{r+v-2}(r+v-2)$. 
Then the homogeneous polynomial $P(\alpha_1,\dotsc,\alpha_k)
\eqdef (-1)^{r+v-1}(P_{r+v}+\alpha_1 P_{r+v-1})$ satisfies $P(1,\dotsc,1)=1$ and vanishes
to order $v$ at $(\alpha_1,\dotsc,\alpha_k)$.

Next we deduce the case of a general $l$ from the case $l=1$. Indeed, suppose $M$ is an $n$-by-$n$ \hbox{$L$-matrix}
over a field $F$ of rank $r$, and $n\geq \binom{r+l-1}{l}+v$. Consider the field $F(X_1,\dotsc,X_l)$, where $X_1,\dotsc,X_l$
are independent indeterminants. As a rank of a matrix does not change when passing to a larger field,
we may treat $M$ as a matrix over $F(X_1,\dotsc,X_l)$. Define a polynomial $g$ by $g(y)=\sum_{i=1}^l X_i y^i$.
By Lemma~\ref{lem:uppbound}, the rank of $g[M]$ does not exceed $\binom{r+l-1}{l}-1$. 

Apply the case $l=1$ to the $g(L)$-matrix $g[M]$ to obtain an integral polynomial $P$. 
Its degree is at most $\binom{r+l-1}{l}+v$. Let $f(x)=\sum_{i=1}^l b_ix^i$ be an arbitrary
polynomial of degree at most $l$ with vanishing constant term. Let $b\eqdef (b_1,\dotsc,b_l)$.
Since $P$ vanishes at $\bigl(g(\alpha_1),\dotsc,g(\alpha_k)\bigr)$
to order at least $v$, it also vanishes at $\bigl(g(\alpha_1),\dotsc,g(\alpha_k)\bigr)|_{X=b}
=\bigl(f(\alpha_1),\dotsc,f(\alpha_k)\bigr)$
to order at least $v$.
\end{proof}

\section{Multivariate polynomials vanishing to high order at a point}\label{sec:vanishing}
A single-variable polynomial of degree $d$ with integer coefficients can vanish at a point $\alpha$ to order exceeding $d/2$ only if
$\alpha\in \Q$. Furthermore, if the polynomial is monic, then $\alpha\in\Z$. 
In this section we prove a generalization of these assertions to homogeneous polynomials in several variables. The following
is the main result of this section.
\begin{lemma}\label{lem:vanishing}
Let $F$ be a field, and let $\overline{F}$ be its algebraic closure. 
Suppose $\alpha=(\alpha_1,\dotsc,\alpha_k)\in \overline{F}^k$ is an arbitrary point, and  $P(x_1,\dotsc,x_k)$ is a homogeneous polynomial with integer coefficients such that
\begin{enumerate}
\item $P$ vanishes at $\alpha$ to an order exceeding $\frac{k-1}{k}\deg P$, and
\item $P(1,\dotsc,1)=1$.
\end{enumerate}
Then there exists a \emph{linear} homogeneous polynomial $Q(x_1,\dotsc,x_k)$ with integer coefficients such that
\begin{enumerate}
\item $Q$ vanishes at $\alpha$, and
\item $Q(1,\dotsc,1)=1$.
\end{enumerate}
\end{lemma}

As a first step, we reformulate the lemma as a result about affine polynomials. So, what we will really prove 
is the following.
\begin{lemma}\label{lem:affinevanishing}
Let $F$ be a field, and let $\overline{F}$ be its algebraic closure.
Suppose $\alpha=(\alpha_1,\dotsc,\alpha_k)\in \overline{F}^k$ is an arbitrary point, and $P(x_1,\dotsc,x_k)$ is a polynomial with integer coefficients such that
\begin{enumerate}
\item $P$ vanishes at $\alpha$ to an order exceeding $\frac{k}{k+1}\deg P$, and
\item $P(1,\dotsc,1)=1$.
\end{enumerate}
Then there exists a degree-one polynomial $Q(x_1,\dotsc,x_k)$ with integer coefficients such that
\begin{enumerate}
\item $Q$ vanishes at $\alpha$, and
\item $Q(1,\dotsc,1)=1$.
\end{enumerate}
\end{lemma}
\begin{proof}[Proof that Lemma~\ref{lem:affinevanishing} implies Lemma~\ref{lem:vanishing}]
Let $P$ be a homogeneous polynomial of degree $d$ satisfying Lemma~\ref{lem:vanishing}.
As the case $\alpha=0$ is trivial, we may assume that $\alpha\neq 0$. Without loss, $\alpha_k\neq 0$.
Define $P'(x_1,\dotsc,x_{k-1})\eqdef P(x_1,\dotsc,x_{k-1},1)$, and $\alpha'=(\alpha_1/\alpha_k,\dotsc,\alpha_{k-1}/\alpha_k)$.
The conclusion of Lemma~\ref{lem:vanishing} then follows from Lemma~\ref{lem:affinevanishing} applied
to $P'$ and $\alpha'$.
\end{proof}

In the case $k=1$, Lemma~\ref{lem:affinevanishing} is a simple consequence of Gauss's lemma. Indeed, 
we may assume that $P$ is primitive, i.e., the coefficients of $P$ are coprime.
Let $P=P_1\dotsc P_l$ be a factorization of $P$ over $\Z$. By Gauss's lemma, the factors $P_1,\dotsc,P_l$ are in fact irreducible over~$\Q$.
Without loss each of $P_1,\dotsc,P_l$ vanishes at $\alpha$. Since $l>\deg P/2$, and $\deg P=\sum \deg P_i$,
at least one of the factors is linear, and the result follows.

For $k\geq 2$, I do not know any equally direct proof. 
The reason is that $Q$ need not be a factor of $P$. Indeed, consider 
the polynomial $x^n+y^n-xy^{n-1}$ and the point $\alpha=(0,0)$.
The polynomial is irreducible because its Minkowski polygon is not a sum
of two smaller lattice polygons \cite[Theorem VI]{ostrowski}.

So we proceed indirectly. We first reduce the lemma to the case $\alpha\in F^k$,
and then characterize those $\alpha\in F^k$ which do not admit polynomial $Q$ as in the lemma. We then show
all these $\alpha$'s do not admit a polynomial $P$. For convenience, we shall treat
the cases $\fchar F=0$ and $\fchar F>0$ separately, the latter case being easier.

\begin{lemma}\label{lem:interpolationone}
Suppose $P(x_1,\dotsc,x_k)$ is a degree-$d$ polynomial which vanishes at points $p_1$ and $p_2$
to orders $d-m_1$ and $d-m_2$ respectively. Then $P$ vanishes at each point of the line $p_1p_2$
to order at least~$d-m_1-m_2$.
\end{lemma}
\begin{proof}
We need to show that, for every multiindex $i$ satisfying $\abs{i}<d-m_1-m_2$, the Hasse derivative
$P^{(i)}(x)$ vanishes at all points of the line $p_1p_2$. By \cite[Lemma~5]{dvir_kopparty_saraf_sudan_multiplicities} the $P^{(i)}$ 
vanishes at $p_1$ and $p_2$ to orders at least $d-m_1-\abs{i}$
and $d-m_2-\abs{i}$ respectively.
The restriction of $P^{(i)}$ onto the line $p_1p_2$ is a univariate polynomial of degree at most $\deg P^{(i)}\leq d-\abs{i}$.
Since the total order of vanishing at $p_1$ and at $p_2$ is
at least $(d-m_1-\abs{i})+(d-m_2-\abs{i})>d-\abs{i}$, it follows that $P^{(i)}$ is identically zero on the line $p_1p_2$.
\end{proof}

By a \emph{flat} we mean an affine subspace (=coset of a vector subspace). An $l$-flat is a flat of dimension~$l$.
For a set $V\subset \overline{F}^k$, let $\Sec_l(V)$ be the union of the all flats spanned by at most $l+1$ points of~$V$.
When $V$ is a variety, then the Zariski closure of $\Sec_l(V)$ is the $l$-th secant variety of~$V$.
By repeatedly applying the preceding lemma we deduce the following.
\begin{corollary}\label{cor:interpolation}
Suppose $P(x_1,\dotsc,x_k)$ is a degree-$d$ polynomial. For an integer $m\geq 0$ let $V_m$ be the set
of points where $P$ vanishes to order exceeding $\frac{m}{m+1}d$. Then $\Sec_l(V_m)\subset V_{m-l}$ for all $l\leq m$.
In particular $P$ vanishes on all points of $\Sec_m(V_m)$.
\end{corollary}

\begin{lemma}\label{lem:subspace}
Suppose $P(x_1,\dotsc,x_k)$ is a polynomial with integer coefficients, and let
$\alpha\in\overline{F}^k$ be a point where $P$ vanishes to order exceeding $\frac{k}{k+1}\deg P$.
Then there exists a flat $V$ defined over $F$ such that
$\alpha\in V$, and
$P$ vanishes on~$V$.
\end{lemma}
\begin{proof}
Let $V_k$ be defined as in Corollary~\ref{cor:interpolation}
The variety $V_k$ is definable over $F$ because it is the intersection
of the zero loci of various Hasse derivatives of~$P$. 
Let $V=\Sec_k(V_k)$. Since $V_k$ is contained in $\overline{F}^k$, it is 
clear that $V$ is the affine span of $V_k$.
By Corollary~\ref{cor:interpolation} $P$ vanishes
on~$V$. As $V_k$ is definable over $F$, then so is $V$.
\end{proof}

We are now ready to prove Lemma~\ref{lem:affinevanishing}. We start with the positive characteristic case.
\begin{proof}[Proof of Lemma~\ref{lem:affinevanishing} in the case $\fchar F=p>0$]
Since polynomial $P$ has integer coefficients and so is defined over every subfield of $F$, we may assume without
loss of generality that $F=\F_p$. Let $V$ be as in Lemma~\ref{lem:subspace}. 
Write $V$ as $V=V_0+v$, where $V_0$ is a vector subspace of $\F_p^k$ and $v\in\F_p^k$. As $P$ vanishes on $V$, but not
at $(1,\dotsc,1)$, it follows that $(1,\dotsc,1)\not\in V$. So, there exists $u\in V_0^{\bot}$
such that $\langle u,(1,\dotsc,1)\rangle\neq \langle u,v\rangle$, and so $x\mapsto \frac{\langle u,x-v\rangle}{\langle u,(1,\dotsc,1)-v\rangle}$
is the desired degree-one polynomial.
\end{proof}

Hence, in the rest of the section we may assume that $\fchar F=0$. In fact, since polynomials $P$ and $Q$
have integer coefficients we may even assume that $F=\Q$.

We shall need the following characterization of when a system of linear
equations in integers admits a solution.
\begin{lemma}[due to van der Waerden, for a proof see \cite{lazebnik_diophantine}]\label{lem:vdW}
Let $M$ be a matrix with rational entries, and $b$ be a rational column vector. Then
the equation $Mz=b$ has an integral vector solution $z$ if and only if, for
every row vector $w^T$ with rational components such that $w^T M$ has integer
components, $w^T b$ is an integer.
\end{lemma}
Note that Lemma~\ref{lem:vdW} was stated in \cite{lazebnik_diophantine} with
an additional restriction that $M$ and $b$ are integral. However, Lemma~\ref{lem:vdW}
follows from the more restrictive version by clearing the denominators.

\begin{lemma}\label{lem:pt}
Suppose $V$ is a flat in $\Q^k$, and suppose that there exists
no degree-one integer polynomial $Q(x_1,\dotsc,x_k)$ vanishing on $V$ and satisfying $Q(1,\dotsc,1)=1$.
Then there exists a point $v\in V\cap \Q^k$ with the property that there exists no degree-one integer polynomial
$Q(x_1,\dotsc,x_k)$ vanishing on $v$ and satisfying $Q(1,\dotsc,1)=1$.
\end{lemma}
\begin{proof}
Write $V$ as $V=\{v_0+t_1 v_1+\dotsb+t_l v_l : t_1,\dotsc,t_l\in \Q\}$, where
$v_0,\dotsc,v_l\in \Q^k$. Treat $v_0,\dotsc,v_l$ as the column vectors, 
and let $\mathbf{1}=(1,\dotsc,1)^T$
denote the all-$1$ column vector. Consider the system of linear equations
in the unknowns $A_0,A_1,\dotsc,A_k$, where we write $A$ to denote the row vector $(A_1,\dotsc,A_k)$,
\begin{align*}
   A_0+A v_0&=0,\\
       A v_i&=0\quad\text{for }i=1,\dotsc,l\\
   A_0+A \mathbf{1}&=1.
\end{align*}
The tuple $(A_0,A_1,\dotsc,A_l)$ is a solution if and only the degree-one polynomial $A_0+A_1x_1+\dotsb+A_kx_k$
vanishes on $V$ and takes value $1$ at $(1,\dotsc,1)$. As we assume that no such polynomial exists, by 
Lemma~\ref{lem:vdW} there exist rational numbers $w_0,w_1,\dotsc,w_l,w'$ such that in the equation
\begin{equation}\label{eq:impossible}
  w_0 A_0 + \sum_{i=0}^l w_i A v_i  +w'(A_0+A\mathbf{1})=w'
\end{equation}
the left side is a linear combination of $A_0,A_1,\dotsc,A_k$ with integer coefficients, whereas
the right side is not an integer. Since the coefficient of $A_0$ is an integer, but $w'$ is not an integer,
it follows that $w_0\neq 0$. Let
\[
  v=v_0+\frac{w_1}{w_0}v_1+\dotsb+\frac{w_k}{w_0}v_k.
\]

We claim that there is no degree-one integer polynomial that vanishes at the point $v$ and takes value $1$
at $(1,\dotsc,1)$. 
Indeed, such a polynomial exists if and only if there is a solution to the system
\begin{align*}
  A_0+A v&=0,\\
  A_0+A \mathbf{1}&=1.
\end{align*}
However, the equation
\[
  w_0(A_0+A v)+w'(A_0+A \mathbf{1})=w'
\]
is the same as \eqref{eq:impossible}, and so the system has no solutions.
\end{proof}

The preceding lemma tells us that for the purpose of proving Lemma~\ref{lem:affinevanishing}
we may assume in effect that $V$ in Lemma~\ref{lem:subspace} is actually a point. The next lemma
deals with that case.

\begin{lemma}\label{lem:ptcond}
Let $v=(\frac{r_1}{s_1},\dotsc,\frac{r_k}{s_k})\in \Q^k$ be a rational point, where $\gcd(r_i,s_i)=1$ for 
all $i=1,\dotsc,k$. There exists a degree-one integer polynomial $Q$ satisfying $Q(v)=0$ and $Q(1,\dotsc,1)=1$
if and only if there exists no prime $p$ such that $r_i\equiv s_i\pmod p$ for all $i=1,\dotsc,k$.
\end{lemma}
\begin{proof}
Consider the equation 
\begin{equation}\label{eq:pt}
  A_1\bigl(1-\frac{r_1}{s_1}\bigr)+\dotsb+A_k\bigl(1-\frac{r_k}{s_k}\bigr)=1\quad\text{ with integer unknowns }A_1,\dotsc,A_k.
\end{equation}
It is easy to see that the solubility of this equation is equivalent to the existence of polynomial~$Q$. If
there is a prime $p$ dividing $s_i-r_i$, then $p\nmid s_i$ since $p\mid s_i,r_i$ contradicts $\gcd(r_i,s_i)=1$.
Hence, $1-r_i/s_i\equiv 0\pmod p$ if we interpret $r_i/s_i$ as a ratio of two elements of $\F_p$.
So, if $p\mid s_i-r_i$ for all $i$, then the equation admits no solution because the left side vanishes modulo~$p$.
Conversely, if $\gcd(s_1-r_1,\dotsc,s_k-r_k)=1$, then there exist integers $B_1,\dotsc,B_k$
such that $B_1(s_1-r_1)+\dotsb+B_k(s_k-r_k)=1$, and so $(A_1,\dotsc,A_k)=(s_1 B_1,\dotsc,s_k B_k)$ is
an integral solution to~\eqref{eq:pt}.
\end{proof}

We are now ready to complete the proof of Lemma~\ref{lem:affinevanishing} in the case $F=\Q$.
Suppose, for the sake of contradiction, that polynomial $P$ satisfies the assumptions of Lemma~\ref{lem:affinevanishing},
but no degree-one polynomial $Q$ fulfilling the conclusion of the lemma exists.
As in the proof of the case $\fchar F>0$, we deduce the existence of a flat $V\subset \{P=0\}\cap \Q^k$
containing $\alpha$.
Lemma~\ref{lem:pt} tells us that there is a point $v\in V$ such that no linear rational polynomial $Q$ satisfying $Q(1,\dotsc,1)$
vanishes at~$v$. Lemma~\ref{lem:ptcond} then yields a prime $p$ such that $v\equiv (1,\dotsc,1)\pmod p$, i.e.,
all the numerators of all the coordinates of $v-(1,\dotsc,1)$ are divisible by $p$.
We conclude that $1=P(1,\dotsc,1)\equiv P(v)=0 \pmod p$, which is a contradiction,
and so Lemma~\ref{lem:affinevanishing} is true after all.

\paragraph{The vanishing condition is optimal}
The order of vanishing in the premise of Lemma~\ref{lem:vanishing} (and hence in Lemma~\ref{lem:affinevanishing}) cannot be reduced.
To see this, we will need a lemma (due to Jacob Tsimerman).
\begin{lemma}
Let $G/F$ be a Galois field extension of degree~$k$. Let $\gamma=(\gamma_1,\dotsc,\gamma_k)\in G^k$ be an arbitrary point.
Let $\gamma^{(1)},\dotsc,\gamma^{(k)}\in G^k$ be the Galois conjugates of $\gamma$.
Then the points $\gamma^{(1)},\dotsc,\gamma^{(k)}$ are linearly independent over $G$ if and only if $\gamma_1,\dotsc,\gamma_k$
are linearly independent over~$F$.
\end{lemma}
\begin{proof}
The `only if' part is trivial, as a linear relation between $\gamma_i$'s also holds between their Galois conjugates.
We shall prove the `if' part.

Let $B_{\gamma}$ be the matrix whose columns are $\gamma^{(1)},\dotsc,\gamma^{(k)}$.
As $G/F$ is Galois, there is an irreducible polynomial $f\in F[x]$ such that $G\cong F[x]/(f)$. Without
loss $G=F[x]/(f)$. Let $\gamma'\eqdef (1,x,\dotsc,x^{k-1})$. The lemma holds for $\gamma'$
because $B_{\gamma'}$ is a Vandermonde matrix.
Since the coordinates of $\gamma$ and $\gamma'$ are $F$-bases for $G$, there is a rational 
invertible matrix $M$ such that $\gamma=M\gamma'$. Since $B_{\gamma'}$ is invertible, then 
so is~$B_{\gamma}=MB_{\gamma'}$.
\end{proof}

Let $G/\Q$ be a Galois extension of degree~$k$. Let $\{1,\gamma_2\dotsc,\gamma_k\}\subset G$
be an integral basis for $G$. Set $\gamma_1\eqdef 1-\sum_{i=2}^k \gamma_i$. Let
$\gamma\eqdef \{\gamma_1,\dotsc,\gamma_k\}$. Let $\gamma=\gamma^{(1)},\dotsc,\gamma^{(k)}$
be the Galois conjugates of $\gamma$. By the previous lemma, $\gamma^{(1)},\dotsc,\gamma^{(k)}$
are linearly independent over~$G$. Define a linear homogeneous polynomial $P_i\in G[x_1,\dotsc,x_k]$ by 
$P_i(x)=\langle x,\gamma^{(i)}\rangle$, and set $P\eqdef \prod_i P_i$.
By the choice of $\gamma$, we have $P(1,\dotsc,1)=1$.
Also note that since $P$ is invariant under $\Gal(G/\Q)$, the coefficients of $P$
are in~$\Q$. Furthermore, since the coefficients of $P_i$ are algebraic integers,
in fact $P\in \Z[x_1,\dotsc,x_k]$.

By the linear independence of $\gamma^{(i)}$'s, the common zero set of $P_1,\dotsc,P_{i-1},P_{i+1},\dotsc,P_k$
is a line through the origin. Let $\alpha^{(i)}$ be any non-zero point on the line.
Note that we may choose $\alpha^{(i)}$'s to be Galois conjugates of one another. The $\alpha^{(i)}$'s are 
linearly independent since $\gamma^{(i)}$'s are.
The polynomial $P$ is of degree $k$ and vanishes at each $\alpha^{(i)}$ to order $k-1$.
However, the conclusion of Lemma~\ref{lem:vanishing} fails for $\alpha^{(1)}$.
Indeed, if there were an integral linear homogeneous polynomial vanishing at $\alpha^{(1)}$,
then it would vanish on all of the $\alpha^{(i)}$'s contrary to the linear independence.

The same construction carries over to finite fields. It is in fact easier as we need not
worry that the coefficients of $\gamma$ are algebraic integers.

\section{Constructions}\label{sec:constr}
In this section we describe the constructions of superlinear-sized matrices for
sets $L$ admitting a primitive linear relation. We then use these to construct
matrices of size $cr$, with $c>1$, for $L$ satisfying a polynomial condition
of Theorem~\ref{thm:minimal}(b), and for sets with $\abs{L}=2$.

\subsection{Construction toolkit}
All our constructions rely on the same basic setup which we describe here.

Let $F$ be a field over which we wish to construct an $L$-matrix.
Let $\F_q$ be a finite field. Let $\Proj^{d-1}(\F_q)$ denote the 
projective space of dimension $d-1$ over $\F_q$. The points of $\Proj^{d-1}(\F_q)$
are the one-dimensional subspaces of $\F_q^d$, and, in general, its $l$-flats 
are $(l+1)$-dimensional subspaces of $\F_q^d$. For a set
$S\subset \F_q^d$ we denote by $\vspan S$ the vector space spanned by~$S$.
When discussing $\Proj^{d-1}(\F_q)$, we shall use concatenation to denote
the span. So, for example if $p,p'$ are two points in $\Proj^{d-1}(\F_q)$,
then $pp'$ is their span, which is a line unless $p=p'$.

Let the Grassmanian $\Gr(s,d)$ be the set of all $s$-dimensional vector subspaces
of $\F_q^d$, or equivalently the set of all $(s-1)$-flats in $\Proj^{d-1}(\F_q)$.
Note that $\Gr(0,d)$ is non-empty, consisting of the unique zero-dimensional subspace of $\F_q^d$;
as an element of $\Proj^{d-1}(\F_q)$ we denote it~$\emptyset$.
Let 
\[
  \Gr(\mathord{\leq} s,d)\eqdef\Gr(0,d)\cup\dotsb\cup \Gr(s,d).
\] 

The following lemma is behind all of our constructions.
\begin{lemma}\label{lem:constr}
Let $\Gr(\mathord{\leq} s,d)$ denote the $\mathord{\leq} s$-dimensional Grassmanian in $\F_q^d$ as defined above.
Let $F$ be a field, $s$ and $d$ be integers satisfying $1\leq s<d$, and suppose $\phi\colon \Gr(\leq s,d) \to F$ is any function.
Then there exists a symmetric $(L,\lambda)$-matrix of size $q^d$ and rank
at most $\abs{\supp \phi}q^s$ with $\lambda=\sum_{W} \phi(W)$
and
\[
   L=\left\{ \sum_{\substack{W\subseteq H\\W\in \Gr(\leq s,d)}} \phi(W) : \text{ hyperplane $H$ in }\Proj^{d-1}(\F_q) \right\},
\]
where the sum is over all flats $W$ of projective dimension less than $s$ (=subspaces $W$ of dimension at most
$s$).  
\end{lemma}
\begin{proof}
While in the application of this lemma it will be easier to use the language of projective geometry,
in the proof of the lemma the language of subspaces will be more convenient.
 
For a subspace $W$ of $\F_q^d$, let $W^{\bot}$ denote the orthogonal complement of $W$.
We will construct an $(L,\lambda)$-matrix whose rows and columns will be
indexed by elements of $\F_q^d$, i.e., a matrix with the underlying vector space  $F^{\F_q^d}$. 
For each subspace $W\in \supp \phi $
and each $y\in \F_q^d$ define a vector $v_y^{(W)}\in F^{\F_q^d}$ by 
\[
  v^{(W)}_{y,x}=
\begin{cases}
\phi(W)&\text{if }x-y\in W^{\bot},\\
0&\text{otherwise}.
\end{cases}
\]
Note that, for a fixed $W$, there are at most $q^{\dim W}$ distinct vectors of the form
$v^{(W)}_y$ as the vector $v_{y}^{(W)}$ depends only on the coset $y+W^{\bot}$.

We then define the matrix $M$ by specifying its rows as
\[
  M_y=\sum_{W \in \supp \phi} v^{(W)}_y.
\]
As its row space is spanned by the vectors of the form $v_y^{(W)}$, the resulting matrix is of rank at most $\sum_W q^{\dim W} \leq \abs{\supp \phi} q^s$.
The diagonal entries are clearly all equal to $\sum_{W} \phi(W)$. More generally, 
the entry in the column indexed by $x$ and the row indexed by $y$ is
\begin{equation}\label{eq:entry}
  \sum_{\substack{W\subseteq (x-y)^{\bot}\\\dim W\leq s}} \phi(W).
\end{equation}
So, the off-diagonal entries belong to the set $L$ defined in the statement of the lemma. From
\eqref{eq:entry} it is also clear that the matrix is symmetric.
\end{proof}

\subsection{Matrices of size \texorpdfstring{$\Omega(r^2)$}{r²} and the case \texorpdfstring{$L\subset \Z$}{L ⊂ Z}}\label{subsec:square}
The next construction is a generalization of the construction of $L$-intersecting
families from \cite{babai_frankl_classification}. Besides being cast in a different setting, 
the version for $L$-intersecting families in \cite{babai_frankl_classification} has an additional requirement that the uniformity 
of the set family is sufficiently large. That is because a matrix $M$ is of rank $r$ if
it factors as $M_1^T M_2$ where $M_1,M_2$ are two $r$-by-$n$ matrices, but $M$
corresponds to a set family only if the entries of $M_1$ and $M_2$ are
nonnegative integers. The nonnegativity constraint is responsible for the extra complexity
in \cite{babai_frankl_classification}.

\begin{theorem}\label{thm:square}
Suppose $L=\{\alpha_1,\dotsc,\alpha_k\}$ is a set admitting a primitive linear relation $A_1\alpha_1+\dotsb+A_k\alpha_k=0$ in which
$k-1$ of the $A_i$'s are nonnegative. Then, for every $r$ there exists a symmetric matrix
of rank $r$ and size $\Omega(r^2)$. In particular, $N(r,L)=\Omega(r^2)$. 
\end{theorem}
\begin{proof}
Without loss $A_2,\dotsc,A_k$ are nonnegative. Let $S=\sum_{i\geq 2} A_i$.
Let $s=1$, $d=2$ and let $q$ be any prime power larger than $S$.
For each $i=2,\dotsc,k$ and each $j=1,\dotsc,A_i$ choose points
$p_{i,j}$ in $\Proj^1(\F_q)$ so that
all these $S$ points are distinct; the choice of $q$ assures that
we can find that many distinct points.
Define the function
$\phi$ by $\phi(p_{i,j})=\alpha_k-\alpha_1$,
and $\phi(\emptyset)=\alpha_1$.
Lemma~\ref{lem:constr} yields a $q^2$-by-$q^2$ matrix $M$ of rank $O(q)$ that is an $(L',\lambda)$-matrix
for
\begin{align*}
  L'&=\alpha_1+\{0,\alpha_2-\alpha_1,\dotsc,\alpha_k-\alpha_1\}=L,\\
  \lambda&=\alpha_1+\sum_{i=2}^k A_i (\alpha_k-\alpha_1)=\sum_{i=1}^k A_i\alpha_i=0.
\end{align*}

We thus obtain a construction of $L$-matrices of size $n$ and rank $O(\sqrt{n})$ whenever $n$ is
a square of a prime power. If $n$ is not a square of a prime power, then we can take
an $n$-by-$n$ submatrix of an $L$-matrix of size $n'$, where $n'$ is the least square of a prime
power satisfying $n'\geq n$. In view of Bertand's postulate, and the fact that matrix rank does not increase by passing to a submatrix,
we obtain a construction for every matrix size.
\end{proof}
\begin{corollary}\label{cor:ltwo}
If $\abs{L}=2$ and $L$ satisfies a primitive linear relation, then $N(r,L)=\Omega(r^2)$.
\end{corollary}
\begin{proof}
In a primitive linear relation $(A_1,A_2)$ of size $2$, one of the $A_1,A_2$ is positive. 
\end{proof}

Because integer vectors satisfy not one, but many linear relations, the preceding theorem implies
a quadratic lower bound for integer sets satisfying a primitive linear relation.
\begin{corollary}
Suppose $L\subset \Z$ is a set satisfying a primitive linear relation, then
$N(r,L)=\Omega(r^2)$.
\end{corollary}
\begin{proof}
Suppose $L=\{\alpha_1,\dotsc,\alpha_k\}$
and $(A_1,\dotsc,A_k)$ is a primitive linear relation, i.e., $\sum A_i=1$ and $\sum A_i\alpha_i=0$.
We may also assume that the relation $(A_1,\dotsc,A_k)$ minimizes the number
of negative coefficients among $A_1,\dotsc,A_k$. If only one of the coefficients
is negative, then the previous theorem applies, and we are done. So, assume, for 
contradiction's sake, that some two coefficients, say $A_1$ and $A_2$ are negative.
Consider the system of linear equations, with unknowns $B_1,B_2,B_3$
\begin{align*}
  0&=B_1+B_2+B_3,\\
  0&=\alpha_1 B_1+\alpha_2 B_2+\alpha_3 B_3.
\end{align*}
It is an underdetermined system of homogeneous equations, and so admits a non-zero solution.
Let $(B_1,B_2,B_3)$ be any solution, which after a suitable scaling we may assume to be integral.
Note that none of $B_1,B_2,B_3$ is zero, for otherwise $\alpha_1,\alpha_2,\alpha_3$ would not be distinct. 
Hence, flipping the signs if necessary, we may also assume that two of the $B_i$'s are positive.
Then the tuple $(A_1+sB_1,A_2+sB_2,A_3+sB_3,A_4,\dotsc,A_k)$ is a primitive linear relation
on $L$, and, for a sufficiently large $s$, two of the first three coefficients are positive.
This contradicts the minimality of $(A_1,\dotsc,A_k)$, implying that only
one of the coefficients is negative after all.
\end{proof}

\subsection{Matrices of size \texorpdfstring{$\Omega(r^{3/2})$}{r**(3/2)} and the case of an arbitrary \texorpdfstring{$L$}{L}}
In the case when $L$ is not a set of integers, we do not have the luxury of choosing
a convenient linear relation, and must make do with a given relation.
\begin{theorem}\label{thm:threehalves}
Let $F$ be a field. Suppose a finite set $L=\{\alpha_1,\dotsc,\alpha_k\}\subset F$ satisfies a primitive linear relation.
Then, for every $r$ there exists a symmetric matrix
of rank at most $r$ and size $\Omega(r^{3/2})$. In particular, $N(r,L)=\Omega(r^{3/2})$.
\end{theorem}
\begin{proof}
Let $\sum A_i \alpha_i=0$, with $A_i\in\Z$ and $\sum A_i=1$, be the primitive linear relation. 
We can rewrite it in the form $\alpha_1+\sum B_{i,i'} (\alpha_i-\alpha_{i'})=0$ where $B_{i,i'}\in \Z_+$.
Let $S=\sum B_{i,i'}$.

We choose $s=2$ and $d=3$, and any $q\geq S-1$. 
Pick a line $l$ in $\Proj^2(\F_q)$.
For each pair $(i,i')$ and for each $j=1,\dotsc,B_{i,i'}$ we shall
choose a distinct point $p_{i,i',j}$ on the line $l$, and a distinct
point $q_{i,i',j}$ not on the line. 
Since $q+1\geq S$, these choices are possible.

Let $l_{i,i',j}$ denote the line spanned by $p_{i,i',j}$ and
$q_{i,i',j}$. We define the non-zero values of the function $\phi$ as follows:
\begin{align*}
  \phi(\emptyset)&=\alpha_1,\\
  \phi(l)&=\sum_{i,i'} B_{i,i'}(\alpha_1-\alpha_{i'}),\\
  \phi(p_{i,i',j})&=\alpha_{i'}-\alpha_1,\\
  \phi(l_{i,i',j})&=\alpha_i-\alpha_{i'}.
\end{align*}
Note that the value of $\phi(l)$ is chosen so that $\phi(\emptyset)+\phi(l)+\sum_{i,i',j} \phi(p_{i,i',j})=\alpha_1$.

We apply Lemma~\ref{lem:constr} to the function $\phi$. We need to verify
that $\sum_{W\subseteq H} \phi(W)\in L$ for every
hyperplane (=line) $H$ in $\Proj^2(\F_q)$.  There are four cases to check.

\textbf{Case 1:} If $H$ contains none of the $p$-points, then $
\sum_{W\subseteq H} \phi(W)=\phi(\emptyset)=\alpha_1$. 

\textbf{Case 2:} If $H$ is the line $l$,
then $\sum_{W\subseteq H} \phi(W)=\phi(\emptyset)+\phi(l)+\sum_{i,i',j} \phi(p_{i,i',j})=\alpha_1$.

\textbf{Case 3:} If $H=l_{i,i',j}$, then $\sum_{W\subseteq H} \phi(W)=\phi(\emptyset)+\phi(p_{i,i',j})+\phi(l_{i,i',j})
=\alpha_i$.

\textbf{Case 4:} If $H$ contains $p_{i,i',j}$, but $H\neq l_{i,i',j}$, then $\sum_{W\subseteq H} \phi(W)=\phi(\emptyset)+\phi(p_{i,i',j})=\alpha_{i'}$.

Finally, we compute the value $\lambda$ in Lemma~\ref{lem:constr} to be
\[
  \lambda=\phi(\emptyset)+\phi(l)+\sum_{i,i',j} \phi(p_{i,i',j})+\phi(l_{i,i',j})=\alpha_1+\sum_{i,i'} B_{i,i'}(\alpha_i-\alpha_{i'})=0.
\]

As in the proof of Theorem~\ref{thm:square}, Bertrand's postulate and rank monotonicity permit us to 
extend the construction from matrices of size $q^3$ to an arbitrary size.
\end{proof}

\subsection{Matrices of size \texorpdfstring{$\Omega(r^{5/3})$}{r**(5/3)} and the case \texorpdfstring{$\abs{L}=3$}{L=3}}
The following is an intermediate result between Theorems~\ref{thm:square} and~\ref{thm:threehalves}. For instance,
it improves upon Theorem~\ref{thm:threehalves} for all sets $L$ of size $\abs{L}\leq 3$.
\begin{theorem}\label{thm:fivethirds}
Let $F$ be a field. Suppose $L=\{\alpha_1,\dotsc,\alpha_k\}\subset F$ is a set admitting a primitive linear relation $A_1\alpha_1+\dotsb+A_k\alpha_k=0$ in which
at least $k-2$ of the $A_i$'s are nonnegative. Then $N(r,L)=\Omega(r^{5/3})$.
\end{theorem}
\begin{proof}
Without loss, $A_3,A_4,\dotsc,A_k>0$. We may also suppose that $A_2<0$. While the case $A_1>0$ is covered by
Theorem~\ref{thm:square}, we make no assumption on $A_1$ as the following proof needs none.
Let $B=-A_2$. Note that $B,A_3,\dotsc,A_k>0$.

We choose $s=3$ and $d=5$ in Lemma~\ref{lem:constr}. We assume that $q$ is large enough to
make the choices described below. Let
$V$ be some $2$-flat in $\Proj^{4}(\F_q)$.
For each $j=1,\dotsc,B$, choose a distinct line $l_j\subset V$ and a  
point $p_j\not\in V$ such that the hyperplanes $Vp_1,Vp_2,\dotsc$ are all distinct. 
Also choose a family of $2$-flats $\{f_{i,j}\}$ (collectively ``$f$-flats'') as follows. 
For each $i=3,4,\dotsc,k$ and for each
$j=1,\dotsc,A_i$, pick a $2$-flat $f_{i,j}$ in $\Proj^{4}(\F_q)$ subject to the three independence conditions:
\begin{enumerate}
\renewcommand*{\theenumi}{I\arabic{enumi}}
\item Any $f$-flat and any $l$-line together span $\Proj^{4}(\F_q)$;
\item No $f$-flat contains any of the $p$-points;
\item Any two $f$-flats span $\Proj^{4}(\F_q)$.
\end{enumerate}

Recall that $l_j p_j$ denotes the $2$-flat spanned by the line $l_j$ and the point $p_j$.
Define the non-zero values of the function~$\phi$ by
\begin{align*}
    \phi(\emptyset)&=\alpha_1,\\
    \phi(l_j)&=\alpha_2-\alpha_1,\\ 
    \phi(V)&=(\alpha_2-\alpha_1)(1-B),\\
    \phi(l_jp_j)&=\alpha_1-\alpha_2,\\
    \phi(f_{i,j})&=\alpha_i-\alpha_1.
\end{align*}
Note that the value of $\phi(V)$ is chosen so that $\phi(\emptyset)+\phi(V)+\sum_j \phi(l_j)=\alpha_2$.

We apply Lemma~\ref{lem:constr}. We need to verify that $\lambda=0$ and 
that $\sum_{W\subseteq H} \phi(W)\in L$ for every
hyperplane $H$ in $\Proj^4(\F_q)$. There are six (easy) cases to check:

\textbf{Case 1:} Suppose $H$ contains some flat $f_{i,j}$. In view of the conditions (I1) and (I3), $H$
contains no $l$-line and no other $f$-flat, respectively. Hence, $\sum_{W\subseteq H} \phi(W)=\phi(\emptyset)+\phi(f_{i,j})=\alpha_i$.

\textbf{Case 2:} Suppose $H$ contains no $f$-flat, and no $l$-line. Then $\sum_{W\subseteq H} \phi(W)=\phi(\emptyset)=\alpha_1$.

\textbf{Case 3:} Suppose $H$ contains two $l$-lines, and contains $p_t$ for some $t$, but no $f$-flat. In view of the condition (I2),
$H$ actually contains all of $V$. Since the hyperplanes $Vp_1,Vp_2,\dotsc$ are all distinct, $H$ contains no $p$-points other than $p_t$. 
Hence, $\sum_{W\subseteq H} \phi(W)=\phi(\emptyset)+\phi(V)+\phi(l_t p_t)+\sum_j \phi(l_j)=\alpha_1$.

\textbf{Case 4:} Suppose $H$ contains two $l$-lines, but no $p$-point or $f$-flat. 
In view of the condition (I2),
$H$ actually contains all of $V$. Hence, $\sum_{W\subseteq H} \phi(W)=\phi(\emptyset)+\phi(V)+\sum_j \phi(l_j)=\alpha_2$.

\textbf{Case 5:} Suppose $l_j\in H$, but $p_j\not\in H$, and $H$ contains no $f$-flat, and no $l$-lines
other than $l_j$. Then
$\sum_{W\subseteq H} \phi(W)=\phi(\emptyset)+\phi(l_j)=\alpha_2$.

\textbf{Case 6:} Suppose $l_j,p_j\in H$, and $H$ contains no $f$-flat, and no $l$-lines other than
$l_j$. Then $\sum_{W\subseteq H} \phi(W)=\phi(\emptyset)+\phi(l_j)+\phi(l_jp_j)=\alpha_1$.\smallskip

Finally, we compute $\lambda$ to be 
\[
  \sum_{W} \phi(W)=\phi(\emptyset)+\phi(V)+\sum \bigl(\phi(l_j)+\phi(l_jp_j)\bigr)+\sum_{i,j}\phi(f_{i,j}) =\alpha_2+(\alpha_1-\alpha_2) B+\sum_i A_i(\alpha_i-\alpha_1)=0.
\]

As in the proof of Theorem~\ref{thm:square}, Bertrand's postulate and rank monotonicity permit us to 
extend the construction from matrices of size $q^5$ to an arbitrary size.
\end{proof}

\subsection{Digraph eigenvalues}\label{subsec:digraph}
We finally have enough tools to construct $\{0,1\}$-matrices with a prescribed eigenvalue of large multiplicity.
We shall not limit ourselves to the set $\{0,1\}$ though, and will present the result in full generality,
for we will also use this construction for the part (c) of Theorem~\ref{thm:minimal}.

\begin{lemma}\label{lem:eigen}
Suppose $L=\{0,\alpha_1,\dotsc,\alpha_k\}$ and let $\widetilde{L}=\{ A_1\alpha_1+\dotsb+A_k\alpha_k : A_1,\dotsc,A_k\in\Z\}$.
Suppose $M$ is an $\widetilde{L}$-matrix of size $n$ with eigenvalue $\lambda$ of multiplicity $m$.
Then for each $l=1,2,\dotsc$ there exists an $L$-matrix $M_l$ of size $ln$ in which $\lambda$ is an eigenvalue of multiplicity~$lm-O(l^{2/3})$.
The constant in the big-oh notation depends on~$L$ and on~$M$.

Furthermore, if $M$ is symmetric, then so is~$M_l$. Also if $\abs{L}=2$, then the exponent $2/3$ can be replaced by~$1/2$.
\end{lemma}
\begin{proof}
Let $I_l$ be the $l$-by-$l$ identity matrix, and put $M_l'=M\otimes I_l$. The multiplicity of $\lambda$ in $M_l'$ is $lm$, and $M_l'$ is of size $ln$.
The matrix $M_l'$ is a block matrix with $n^2$ blocks, each of which is of the form~$\beta I_l$ for various~$\beta\in \widetilde{L}$.

Let $\beta=A_1\alpha_1+\dotsb+A_k\alpha_k$ be an arbitrary element of~$\widetilde{L}$. Then
$(1+\sum A_i,-A_1,\dotsc,-A_k)$ is a primitive linear relation on $\{\beta,\beta+\alpha_1,\dotsc,\beta+\alpha_k\}$,
for $\left(1+\sum A_i\right)\beta-\sum_i A_i(\beta+\alpha_i)=0$. Hence, by Theorem~\ref{thm:threehalves},
there exists a symmetric $\{\beta,\beta+\alpha_1,\dotsc,\beta+\alpha_k\}$-matrix 
$Q_\beta'$ of size $l$ and rank $O(l^{2/3})$; in the case~$\abs{L}=2$, Corollary~\ref{cor:ltwo} guarantees a better bound of~$O(l^{1/2})$.
Let $Q_{\beta}=Q_{\beta}'-\beta J_l$, where $J$ is the $l$-by-$l$ all-$1$ matrix.
As $\rank J_l=1$, the rank of $Q_{\beta}$ is also $O(l^{2/3})$ (resp.~$O(l^{1/2})$).
Note that $Q_{\beta}$ is an $(L,-\beta)$-matrix, and $\beta I_l+Q_{\beta}$ is an $L$-matrix.

We replace each block in $M_l'$ of the form $\beta I_l$ by $\beta I_l+Q_{\beta}$ to obtain matrix $M_l$. 
Each such replacement adds to $M_l'$ a matrix of the same rank as $Q_{\beta}$, namely
the matrix that is all $0$ except for a single block that is~$Q_{\beta}$. Hence,
a single replacement changes the multiplicity of eigenvalue $\lambda$ by at most $O(l^{2/3})$ (resp.~$O(l^{1/2})$).
Since the number of blocks, $n^2$, is constant, the requisite bound on the rank of $M_l$ follows. 

If $M$ is symmetric, then the block structure in $M_l'$ is symmetric, which in view of the $Q_{\beta}$'s being
symmetric implies that the final matrix $M_l$ is symmetric, too.
\end{proof}

Theorem~\ref{thm:quadratic} about the maximum multiplicity of an eigenvalue in a $\{0,1\}$-matrix
is just a simple corollary of the preceding construction.
\begin{proof}[Proof of Theorem~\ref{thm:quadratic}]Let $F_0$ be the prime subfield of~$F$.

\textbf{Part (a)}: If $\alpha$ is in $F_0$, and $1/(1-\alpha)$ is an algebraic integer,
then $1/(1-\alpha)=A$ for some $A\in\Z$. Hence, $(A-1)\cdot 1-A\cdot \alpha=0$ is a primitive relation.
So, for every $r$, an $\{1,\alpha\}$-matrix of size $\Theta(r^2)$ and rank $r$ exists by Theorem~\ref{thm:square}. A quadratic
upper bound follows from Proposition~\ref{prop:basic}. Hence, $N(r,\{1,\alpha\})=\Theta(r^2)$,
of which $E(n,\lambda)=n-\Theta(\sqrt{n})$ is a trivial reformulation, as shown by relations~\eqref{eq:en}.

\textbf{Part (b)}: Suppose $\lambda$ is an algebraic integer with minimal polynomial $f(x)=x^d+\sum_{i=0}^{d-1} a_ix^i=0$ with $a_i\in\Z$ and $d\geq 2$.
Let $M$ be the companion matrix of $f$, which is a $d$-by-$d$ integer matrix whose characteristic polynomial
is~$f$. In particular, $\lambda$ is an eigenvalue of $M$. So, by Lemma~\ref{lem:eigen}, we have $E(ld,\lambda)\geq l-O(l^{1/2})$.
Since  $E(n,\lambda)$ is nondecreasing in $n$, we conclude that
$E(n,\lambda)\geq n/d-O(\sqrt{n})$. In view of \eqref{eq:eigenupper}, and relations~\eqref{eq:en}, the proof is complete.
\end{proof}

\subsection{Graph eigenvalues}\label{subsec:graph}
Recall that we call an algebraic integer $\lambda$ \emph{representable} if there exists an integral symmetric matrix~$M$
whose only eigenvalues are $\lambda$ and its conjugates. Equivalently, the multiplicity of the eigenvalue $\lambda$ in $M$
is equal to exactly $n/d$, where $d$ is the degree of $\lambda$.
When the progenitor matrix $M$ is symmetric, Lemma~\ref{lem:eigen} yields symmetric matrices, and
so $E_s(n,\lambda)\geq n/d-O(\sqrt{n})$ for a representable~$\lambda$. Thus Theorem~\ref{thm:reprconj}
is just a special case of Lemma~\ref{lem:eigen}. 

For all totally real $\lambda$ of degree $d$, Mario Kummer \cite{kummer_symmetric} constructed symmetric integral matrices of size at most $9d$ with eigenvalue $\lambda$. In view
of Lemma~\ref{lem:eigen} this implies that $E_s(n,\lambda)\geq n/9d-O(\sqrt{n})$ for such $\lambda$.

\subsection{Linear-sized matrices from polynomial relations}
In this subsection we give a construction used in Theorem~\ref{thm:minimal}. Namely, we shall show 
that a single polynomial relation on $L$ implies that $N(r,L)\geq cr$ for some $c>1$.
\begin{theorem}\label{thm:polytobetter}
Suppose $L=\{\alpha_1,\dotsc,\alpha_k\}$ and $P$ is a homogeneous polynomial with integer coefficients satisfying
$P(1,\dotsc,1)=1$ and $P(\alpha_1,\dotsc,\alpha_k)=0$. Let $d=\deg P\geq 2$ be the degree of the polynomial. Then $N(r,L)/r\geq 1+1\!\bigm/\!\left[\binom{d+k-2}{k-1}-1\right]+O(r^{-1/3})$ as $r\to\infty$.
\end{theorem}
\begin{proof}
Let $\alpha_i'=\alpha_i-\alpha_1$, and $L'=\{0,\alpha_2',\dotsc,\alpha_k'\}$.
It suffices to construct, for all large $l$, an $L'$-matrix $M$ of size $l\binom{d+k-2}{k-1}$ and eigenvalue $\alpha_1$ of multiplicity $l-O(l^{2/3})$,
for then $M+\alpha_1(J-I)$ is an $L$-matrix of rank at most $l\left[\binom{d+k-2}{k-1}-1\right]+O(l^{2/3})$. 
Let $Q(x_1,\dotsc,x_k)\eqdef P(x_1,x_2+x_1,\dotsc,x_k+x_1)$.
Let $\widetilde{L}=\{A_2\alpha_2'+\dotsb+A_{k}\alpha_{k}' : A_2,\dotsc,A_{k}\in\Z\}$.
In view of Lemma~\ref{lem:eigen}, it suffices to construct an $\widetilde{L}$-matrix of size $\binom{d+k-2}{k-1}$ 
having an eigenvalue~$\alpha_1$.

Let $\mathcal{M}_{d-1},\mathcal{M}_d$ be the families of all homogeneous monomials in $x_1,\dotsc,x_k$ of degrees $d-1$ and $d$ respectively.
We call monomial $m'$ a \emph{predecessor} of a monomial $m$ if $m=x_im'$ for some $i\geq 2$. 
Each monomial $m$ not of the form $x_1^i$ has a predecessor (possibly several). For such an $m$, choose
any predecessor, and denote it by $\pred (m)$. Define the index $i_m$ by $m=x_{i_m} \pred (m)$.
In particular, 
\begin{equation}\label{eq:predalpha}
  m(\alpha_1,\alpha_2',\dotsc,\alpha_k')=\alpha_{i_m}' \pred(m)(\alpha_1,\alpha_2',\dotsc,\alpha_k').
\end{equation}

Let $Q(x_1,\dotsc,x_k)=\sum_m c_m m$ be the expansion of $Q$ as a linear combination of monomials in~$\mathcal{M}_d$.  We shall construct an $\widetilde{L}$-matrix $M$ whose rows and columns are indexed by $\mathcal{M}_{d-1}$. 
The matrix $M$ is a sum of two matrices $M=M^{(1)}+M^{(2)}$, whose non-zero entries are defined to be
\begin{align*}
  M^{(1)}_{\pred (m' x_1),m'} &=\alpha_{i_{m' x_1}}'&&\text{ for each }m'\in \mathcal{M}_{d-1}\setminus \{x_1^{d-1}\},\\
  M^{(2)}_{m,x_1^{d-1}}&=\!\!\!\sum_{\substack{\bar{m}\in \mathcal{M}_d\setminus\{x_1^d\}\\\pred(\bar{m})=m}}\!\!(-c_{\bar{m}} \alpha_{i_{\bar{m}}}')&&\text{ for each }m\in \mathcal{M}_{d-1}\setminus\{x_1^{d-1}\}.
\end{align*}
The matrix $M$ has $\alpha_1$ as an eigenvalue. Indeed, letting $M^{(3)}=M^{(1)}-\alpha_1 I$ we obtain from~\eqref{eq:predalpha}
\[
 \sum_{m\in\mathcal{M}_{d-1}} m(\alpha_1,\alpha_2',\dotsc,\alpha_k') M^{(3)}_{m,m'}=\begin{cases}-\alpha_1^d&\text{if }m'=x_1^{d-1},\\0&\text{otherwise.}\end{cases}
\]
Similarly from \eqref{eq:predalpha} we deduce that  $\sum m(\alpha_1,\alpha_2',\dotsc,\alpha_k') M^{(2)}_{m,x_1^{d-1}}=\alpha_1^d-P(\alpha_1,\alpha_2',\dotsc,\alpha_{k-1}')$.
It then follows that $M-\alpha_1 I=M^{(1)}+M^{(3)}$ is singular, and so $M$ has eigenvalue $\alpha_1$.
\end{proof}

\section{Proofs of Theorems~\ref{thm:minimal}, \ref{thm:superlinear} and \ref{thm:lineareq}}
\label{sec:main}
In this section we reap the fruits of the work above, and derive the main results of this paper.

For a matrix $M$ let $(M\ 1)$ denote the matrix obtained from $M$ by appending an all-$1$ column.
Let 
\[
  N_0(r,L)=\max\{n : M\text{ is an }n\text{-by-}n\ L\text{-matrix with }\rank(M\ 1)\leq r\}.
\]
In the proof of Theorem~\ref{thm:minimal} we will need the following easy fact.
\begin{lemma}\label{lem:superadd}
Let $L$ be an arbitrary finite subset in some field. Then
\[
  N_0(r_1+r_2,L)\geq N_0(r_1,L)+N_0(r_2,L).
\]
\end{lemma}
\begin{proof}
Let $n_1=N_0(r_1,L)$ and $n_2=N_0(r_2,L)$.
Let $M_1$ and $M_2$ be square matrices
of dimensions $n_1$ and $n_2$ satisfying $\rank(M_i\ 1)\leq r_i$ for $i=1,2$.
Pick any $\alpha\in L$, and let $M$ be the block matrix $\left(\begin{smallmatrix}M_1&\alpha J\\\alpha J&M_2\end{smallmatrix}\right)$.
Every linear relation satisfied by the rows of $(M_1\ 1)$ is satisfied
by the first $n_1$ rows of $(M\ 1)$. Similarly, every linear relation satisfied by the rows
of $(M_2\ 1)$ is satisfied by the last $n_2$ rows of $(M\ 1)$.
Hence, the matrix $M$ is a witness to $N_0(r_1+r_2,L)\geq n_1+n_2$.
\end{proof}

We are now ready to prove Theorem~\ref{thm:minimal} characterizing those sets $L$ for which $N(r,L)$ is $r+O(1)$.
\begin{proof}[Proof of Theorem~\ref{thm:minimal}]\ 

\makebox[5.5em][l]{\textbf{(a)$\implies$(b)}:} This is the special case $(l,v)=(1,1)$ of Theorem~\ref{thm:genupper}.

\makebox[5.5em][l]{\textbf{(b)$\implies$(c)}:} If the polynomial $P$ is linear, then this follows from Theorem~\ref{thm:threehalves}.
If the polynomial $P$ has degree $d\geq 2$, then this is the content of Theorem~\ref{thm:polytobetter}.

\makebox[5.5em][l]{\textbf{(c)$\implies$(a)}:} This is trivial.\medskip

Since $N(r,L)\geq N_0(r,L)$ and $N_0(r-1,L)\geq N(r,L)$,
the limit of $N(r,L)/r$ as $r\to\infty$ is equal to the limit of $N_0(r,L)/r$. The latter exists as a consequence of superadditivity of $N_0(r,L)$ (Lemma~\ref{lem:superadd}).
\end{proof}

\begin{proof}[Proof of Theorem~\ref{thm:superlinear}]\

\makebox[5.5em][l]{\textbf{(a)$\implies$(b)}:} If $N(r-1,L)\geq kr+1$ for some $r$, then by Theorem~\ref{thm:genupper}
with $(l,v)=\bigl(1,(k-1)r+1\bigr)$, we see that there is a homogeneous polynomial $P$ of degree at most~$kr+1$
vanishing to order at least $(k-1)r+1$ at $(\alpha_1,\dotsc,\alpha_k)$ and satisfying $P(1,\dotsc,1)=1$.
By Lemma~\ref{lem:vanishing}, $L$ in fact satisfies a primitive linear relation.

\makebox[5.5em][l]{\textbf{(b)$\implies$(c)}:} This is the content of Theorem~\ref{thm:threehalves}. The `furthermore' part
is the content of Theorems~\ref{thm:square} and~\ref{thm:fivethirds}.

\makebox[5.5em][l]{\textbf{(c)$\implies$(a)}:} This is trivial.
\end{proof}

We next tackle Theorem~\ref{thm:lineareq}, asserting that the growth of $N(r,L)$ is determined 
by the primitive linear relations on~$L$. For $L=\{\alpha_1,\dotsc,\alpha_k\}$ we let $R(L)\eqdef\nobreak \{(A_1,\dotsc,A_k) : A_1\alpha_1+\dotsb+ A_k\alpha_k=0\}$ (resp.~$P(L)\eqdef\{ (A_1,\dotsc,A_k)\in R(L) : A_1+\dotsb+A_k=1\}$) to be the collection of all (resp.~all primitive) linear relations on $L$.
\begin{lemma}\label{lem:reli}
If $L$ and $L'$ are two sets such that $P(L)=P(L')$, then either $P(L)=P(L')=\emptyset$, or
$R(L)=R(L)$.
\end{lemma}
\begin{proof}
Suppose that $P(L)=P(L')\neq \emptyset$.  Let $B=(B_1,\dotsc,B_k)\in P(L)$ be any primitive relation,
and suppose $C=(C_1,\dotsc,C_k)\in R(L)$. Then $C+tB\in R(L)$ for every $t\in\Z$,
and in particular for $t=1-(C_1+\dotsb+C_k)$, in which case $C+tB\in P(L)$. Hence $C+tB\in P(L')$ and thus
$C=(C+tB)-tB\in R(L')$.
\end{proof}
\begin{lemma}\label{lem:gamma}
Let $F_{\text{big}}/F_{\text{small}}$ be a finite field extension of degree $D$, and assume that
$\{\beta_1,\dotsc,\beta_d\}\subset\nobreak F_{\text{big}}$ is a non-empty set that is linearly independent 
over~$F_{\text{small}}$. Then there exists
a set $\{\gamma_1,\dotsc,\gamma_s\}$ of size $s\geq D/2d$
such that the $sd$ products $\{\gamma_i \beta_j : i=1,\dotsc,s\ j=1,\dotsc,d\}$ are linearly
independent over~$F_{\text{small}}$.
\end{lemma}
\begin{proof}
Let $\{\gamma_1,\dotsc,\gamma_s\}$ be a maximal set satisfying the conclusion of the lemma.
By maximality for every $\gamma\in F_{\text{big}}$, we have a relation
of the form $\sum_j c_j \beta_j \gamma=\sum_{i,j} c_{i,j} \gamma_i \beta_j$ for some
$c_j,c_{i,j}\in F_{\text{small}}$ and with not all $c_j$ being zero. Hence, $F_{\text{big}}$ is equal to 
\[ 
  I\eqdef\left\{ \frac{\sum_{i,j} c_{i,j} \gamma_i \beta_j}{\sum_j c_j \beta_j} : c_j, c_{i,j}\in F_{\text{small}} \right\}.
\]
We claim that $d(s+1)-1\geq D$. It is easiest to see this if $F_{\text{small}}$
and $F_{\text{big}}$ are finite fields. In that case, the cardinality of $I$ is at most
$(\abs{\F_{\small}}^{d(s+1)}-1)/(\abs{F_{\text{small}}}-1)$, where the 
term $\abs{\F_{\small}}^{d(s+1)}-1$ counts the number of ways to choose $c$'s so that not all
of them are zero, and the factor of $1/(\abs{F_{\text{small}}}-1)$
is due to the homogeneity in $c$'s. The claim then follows from $\abs{I}=\abs{F_{\text{big}}}=\abs{F_{\text{small}}}^{D}$.

If both $F_{\text{small}}$ and $F_{\text{big}}$ are infinite, we
identify $F_{\text{big}}\otimes_{F_{\text{small}}} \overline{F_{\text{small}}}$ with an affine space of dimension $[F_{\text{big}}:F_{\text{small}}]$
over~$F_{\text{small}}$. Under this identification the set $I$ becomes a set
of $F_{\text{small}}$-points of a variety of dimension at most~$d(s+1)-1$. If $S\subset F_{\text{small}}$ is any set
of size $N$, then the Schwartz--Zippel lemma for varieties \cite[Lemma~14]{bukh_tsimerman} tells us that
$\abs{S^D\cap I}=O(N^{d(s+1)-1})$. Since $S^D\cap I=S^D\cap F_{\text{big}}=S^D$ and $F_{\text{small}}$ contains arbitrarily
large sets, the claim follows (in the infinite field case).

In either case, from $d(s+1)-1\geq D$ we deduce that $ s\geq \lfloor D/d \rfloor\geq D/2d$.
\end{proof}
\begin{lemma}\label{lem:fieldchange}
Suppose $F_{\text{big}}/F_{\text{small}}$ is a finite field extension, and $v_1,\dotsc,v_r$
are vectors in $F_{\text{big}}^n$. Suppose the components of $v_1,\dotsc,v_r$
span a vector space of dimension $d$ over~$F_{\text{small}}$. Let
$V_{\text{small}}$ and $V_{\text{big}}$ be the spans of $v_1,\dotsc,v_r$
over $F_{\text{small}}$ and over $F_{\text{big}}$ respectively.
Then $\dim_{F_{\text{small}}} V_{\text{small}}\leq 2d \dim_{F_{\text{big}}} V_{\text{big}}$.
\end{lemma}
\begin{proof}
Let $\{\beta_1,\dotsc,\beta_d\}$ be a basis for the vector space spanned
by the components of $v_1,\dotsc,v_r$ over~$F_{\text{small}}$.
Then $V_{\text{small}}\subset \bigoplus_i \beta_i F_{\text{small}}^n$. Let $\{\gamma_1,\dotsc,\gamma_s\}$
be as in Lemma~\ref{lem:gamma}. Then the vector spaces $\gamma_1 V_{\text{small}},\dotsc,\gamma_s V_{\text{small}}$
are linearly independent over $F_{\text{small}}$. As they are subspaces of $F_{\text{big}}$, we infer that
\[ 
  s\dim V_{\text{small}}\leq \dim_{F_{\text{small}}} V_{\text{big}}=[F_{\text{big}}:F_{\text{small}}]\dim_{F_{\text{big}}} V_{\text{big}}.
\]
As $s\geq [F_{\text{big}}:F_{\text{small}}]/2d$, the lemma follows.
\end{proof}
\begin{proof}[Proof of Theorem~\ref{thm:lineareq}]
By the assumption $P(L)=P(L')$. If in addition, $P(L)=P(L')=\emptyset$, then by Theorem~\ref{thm:superlinear}, $r\leq N(r,L),N(r,L')\leq kr+k+1\leq 2kr$.
So, assume that $P(L)=P(L')\neq\emptyset$. By Lemma~\ref{lem:reli}, we then conclude that $R(L)=R(L')$.

By rescaling $L$ and $L'$ as necessary, we may assume that $1\in L,L'$. Rescaling changes neither $N(r,L),N(r,L')$
nor $R(L),R(L')$.

Consider set $L$, and inside $L$ consider a maximal subset that satisfies no integer relation. By relabeling
elements of $L$ if necessary, we may assume that the subset is $\{1,\alpha_2,\dotsc,\alpha_l\}$ and $\alpha_{l+1},\dotsc,\alpha_k$
are the remaining elements of~$L$. By the maximality assumption, there exist rational linear forms $f_{l+1},\dotsc,f_k$ such that
$\alpha_i=f_i(1,\alpha_2,\dotsc,\alpha_l)$ for $i=l+1,\dotsc,k$. Note that linear relations
$\alpha_i-f_i(1,\alpha_2,\dotsc,\alpha_l)=0$ consistute a basis for the $\Q$-vector space of all linear relations
among $\alpha_2,\dotsc,\alpha_k$.

Let $F$ be the field containing $L$ and $L'$, and let $F_0$ be the prime subfield of $F$.
Let $\widetilde{F}\eqdef F_0(x_2,\dotsc,x_l)$. For $i=l+1,\dotsc,k$ put $x_i\eqdef f_i(1,x_2,\dotsc,x_l)$, and
let $\widetilde{L}\eqdef \{1,x_2,\dotsc,x_k\}$. Note that $R(\widetilde{L})=R(L)$ because
relations $x_i-f_i(1,x_2,\dotsc,x_l)$ consistute a basis for the $\Q$-vector space of all linear
relations among $x_2,\dotsc,x_k$. 

\textsc{Claim 1:} $N(r,\widetilde{L}_{\widetilde{F}})\leq N(r,L)$.

\textsc{Claim 2:} $N(r,L)\leq N(2tr,\widetilde{L}_{\widetilde{F}})$ for some $t\leq l$.

To complete the proof it suffice to prove these two claims. Indeed,
as the field $\widetilde{F}$ and the set $\widetilde{L}$ depend only on $R(L)$,
and $R(L)=R(L')$, if the inequalities in Claims 1 and 2 hold for $L$, they also hold for $L'$.
The claims then imply $N(r,L)\leq N(2tr,\widetilde{L}_{\widetilde{F}})\leq N(2tr,L')$, and similarly
with the roles of $L$ and $L'$ swapped.\smallskip

\textsc{Proof of Claim 1:} Given an $\widetilde{L}$-matrix $\widetilde{M}$, we can define an $L$-matrix
$M$ by replacing each entry $x_i$ in $\widetilde{M}$ by $\alpha_i$. Since every linear
relation satisfied by the rows of $\widetilde{M}$ is also satisfied by the corresponding
rows of $M$, it follows that $\rank M\leq \rank \widetilde{M}$.\smallskip

\textsc{Proof of Claim 2:} Consider the maximal subset of $\{1,\alpha_2,\dotsc,\alpha_l\}$
that is algebraically independent. Without loss of generality, it is $\{1,\alpha_{t+1},\dotsc,\alpha_l\}$
for some~$t$. The field $F_{\text{big}}\eqdef \F_0(\alpha_2,\dotsc,\alpha_l)$ is a finite algebraic
extension of $F_{\text{small}}\eqdef F_0(\alpha_{t+1},\dotsc,\alpha_l)$. 

Suppose $M$ is an $L$-matrix, and let $\widetilde{M}$ be the $\widetilde{L}$-matrix obtained from $M$
by replacing each entry~$\alpha_i$ by~$x_i$.  Let~$r$ be the rank of $M$ over $F$. Note that $r$ is also
the rank of $M$ over $F_{\text{big}}$. Let $V_{\text{big}}$ and $V_{\text{small}}$ be the spans of the rows 
of $M$ over $F_{\text{big}}$ and $F_{\text{small}}$ respectively. Since the entries of $M$
are spanned by $\{1,\alpha_2,\dotsc,\alpha_t\}$ over $F_{\text{small}}$, 
from Lemma~\ref{lem:fieldchange} we deduce that $\dim V_{\text{small}}\leq 2tr$.
As $F_{\text{small}}$ is naturally isomorphic to $F_0(x_{t+1},\dotsc,x_l)$, 
any linear relation between rows of $M$ with coefficients in $F_{\text{small}}$ corresponds
to a linear relation between rows of $\widetilde{M}$. Hence, $\rank \widetilde{M}\leq 2t\rank M$.
\end{proof}

\section{Remarks and open problems}
\begin{itemize}
\item
I know only one example of a $k$-element set $L$ that attains the bound
$N(r,L)\leq r^k/k!+O(r^{k-1})$ of Proposition~\ref{prop:basic} without the loss of a multiplicative constant. That
set is $L=\{1,2,\dotsc,k\}$ and its multiples. Namely, let $A$ be the $r$-by-$\binom{r}{k}$ 
matrix whose columns are the characteristic vectors of the $k$-elements subsets of 
a fixed $r$-element set. Then $kJ-A^TA$ is an $L$-matrix of dimension $\binom{r}{k}$ with $L=\{1,2,\dotsc,k\}$.
Its rank is at most~$r+1$.

It would be very interesting to decide if there are any other examples that attain the bound in Proposition~\ref{prop:basic}.

\item For each $l$ and $r$ there exist a $k$-uniform $\{0,1\}$-intersecting family $\mathcal{F}$ of subsets of $[r]$ with
$\abs{\mathcal{F}}\geq \binom{r}{2}/\binom{l}{2}+O(n)$ (see \cite{martin_rodl} for a particularly simple construction).
That implies the bound $N(r,\{l-1,l\})\geq \binom{r}{2}/\binom{l}{2}+O(r)$. Can this be improved?

Interestingly, almost the same bound, namely $N(r,\{l-1,l\})\geq (r/l+1)^2+O(r)$ can be obtained very differently. Namely,
one can use the relation between graph eigenvalues of $N(r,L)$ for two element sets $L$ in \eqref{eq:en}.
To get the stated bound one uses the square lattice graphs, which are strongly regular graphs with 
parameters $\bigl(n^2,l(n-1),(l-1)(l-2)+n-2,l(l-1)\bigr)$; see \cite{cameron_strongly_regular_survey} for a definition of
these graphs and a survey of strongly regular graphs in general.

\item In this paper we focused on the magnitude of the leading term in the asymptotics for $N(r,L)$.
However, in applications even the lower-order terms in \eqref{eq:basicupperbound} are of much interest, see 
for example \cite[Theorem 4.1.1]{blokhuis_fewdist} or \cite[Theorem~1]{musin_twodistance} (and its generalization in
\cite[Theorem~3.2]{nozaki_shinohara_sdistance}). I do not know if lower-order improvements to bounds in this paper are 
possible.

In particular, is it possible to show, at least for some $\lambda$, that the maximum multiplicity of an eigenvalue $\lambda$ 
in an $n$-vertex graph satisfies $E_s(n,\lambda)\leq \frac{n}{\deg \lambda}-c\sqrt{n}$ for some $c>0$ and all $n\geq n_0(\lambda)$?

\item I conjecture that the exponent $3/2$ in Theorem~\ref{thm:threehalves} can be replaced by $2$. Namely, any set $L$
admitting a primitive linear relation satisfies $N(r,L)=\Omega(r^2)$. 

As evidence, here is a construction showing that $N(r,\{x+y,3x,3y\})=\Omega(r^2)$ for any $x,y\in F$.
Note that $3\cdot(x+y)-1\cdot 3x-1\cdot 3y=0$ and the relation $(3,-1,-1)$ is not covered by Theorem~\ref{thm:square}.
Let $p_1,\dotsc,p_4$ any four points in $\Proj^{3}(\F_q)$ that span $\Proj^{3}(\F_q)$, and define
the function $\phi \colon \Gr(\leq 2,4)\to F$ by
\begin{align*}
  \phi(\emptyset)&=x+y,\\
  \phi(p_1)=\phi(p_3)&=2x-y,\\
  \phi(p_2)=\phi(p_4)&=2y-x,\\
  \phi(p_1p_2)=\phi(p_2p_4)=\phi(p_3p_4)&=x-2y,\\
  \phi(p_2p_3)=\phi(p_1p_3)=\phi(p_1p_4)&=y-2x.
\end{align*}
Lemma~\ref{lem:constr} applied to this $\phi$ shows that $N(r,\{x+y,3x,3y\})=\Omega(r^2)$.

\end{itemize}

\textbf{Acknowledgements.} I am thankful
to James Cummings for discussions, encouragement, and for careful reading of an earlier version of this paper.
I am grateful to MathOverflow users David Lampert, Will Sawin 
and Dracula for help \cite{sawin,dracula} with Lemma~\ref{lem:vanishing}. I thank Gary Greaves
for bringing to my attention the references \cite{estes_guralnick,dobrowolski_symmetric,mckee_six}. I also benefited from
discussions with Tibor Szab\'o, Jacob Tsimerman and Michel Waldschmidt, and from comments of the anonymous referee.

\bibliographystyle{plain}
\bibliography{lmatrices}

\end{document}